\documentclass{article}

\usepackage[english]{babel}

\usepackage[english]{babel}

\usepackage[a4paper,top=2cm,bottom=2cm,left=3cm,right=3cm,marginparwidth=1.75cm]{geometry}

\usepackage{graphicx,amssymb,lineno}
\usepackage{amsmath}
 \allowdisplaybreaks
\usepackage[colorlinks=true, allcolors=blue]{hyperref}
\newtheorem{rmk}{Remark}
\newtheorem{theorem}{Theorem}

\begin{document}

\title{A new thin layer model for viscous flow between two nearby non-static surfaces}
\author{Jos\'e M. Rodr\'{\i}guez${}^{\ref{dir1}, \ref{dir2}}$, Raquel Taboada-Vázquez${}^{\ref{dir1}, \ref{dir3}}$}
\date{}

\maketitle

{\footnotesize 
\begin{enumerate}
\item CITMAga, Santiago de Compostela, Spain. \label{dir1}
\item Department of Mathematics, Higher Technical University College of Architecture, Universidade da Coruña, A Coruña, Spain. E-mail: jose.rodriguez.seijo@udc.es. \label{dir2}
\item Department of Mathematics, School of Civil Engineering, Universidade da Coruña, A Coruña, Spain. E-mail: raquel.taboada@udc.es. \label{dir3}
\end{enumerate}
}

\begin{abstract}
We propose a two-dimensional flow model of a viscous fluid between two close moving surfaces. We show, using a formal asymptotic expansion of the solution, that its asymptotic behavior, when the distance between the two surfaces tends to zero, is the same as that of the the Navier-Stokes equations.

The leading term of the formal asymptotic expansions of the solutions to the new model and Navier-Stokes equations are solution of the same limit problem, and the type of the limit problem depends on the boundary conditions. If slip velocity boundary conditions are imposed on the upper and lower bound surfaces, the limit is a solution of a lubrication model, but if the tractions and friction forces are known on both bound surfaces, the limit is a solution of a thin fluid layer model.

The model proposed has been obtained to be
 a valuable tool for computing viscous fluid flow between two nearby moving surfaces, without the need to decide a priori whether the flow is typical of a lubrication or a thin fluid layer problem, and without the enormous computational effort that would be required to solve the Navier-Stokes equations in such a thin domain.
\end{abstract}


\section{Introduction}\label{sec1}

In our previous work \cite{RodTabJMAA2021}, we used the asymptotic expansions technique to study the behavior of a viscous fluid that flows between two very close moving surfaces. Asymptotic analysis is a mathematical tool that has been used successfully (since the pioneering works of Dean \cite{Dean1}-\cite{Dean2},  Friedrichs and Dressler \cite{FriedrichsDressler} and Goldenveizer \cite{Goldenveizer}) to obtain and justify mathematical models, in solid mechanics \cite{Rigolot1972}-\cite{TutekAganovicNedelec} and fluid mechanics \cite{Cimatti1983}-\cite{SG3}, when at least one of the dimensions of the domain is much smaller than the others. Using the same mathematical technique, the authors have also proposed several new shallow water models \cite{RTV1}-\cite{RTV6} and  curved-pipe flow models \cite{CR1}-\cite{CR2}.

We observed, in our prior article \cite{RodTabJMAA2021},  that the viscous fluid that moves between two nearby surfaces has two very different behaviors, depending on the boundary conditions of the problem. If the pressure differences are large in the open part of the domain boundary (that is, the region of the domain boundary between the two surfaces), then the fluid obeys equation \eqref{Reynolds_gen}, which resembles a lubrication problem. If the pressure differences are small in the mentioned region of the domain boundary, then the fluid obeys equation \eqref{ec_Vi0-v2}, which is a thin fluid layer problem (it can also be  understood as a shallow water problem in which the depth is known).

This behavior reminds us of that observed in the works of Ciarlet et al. \cite{CiarletLodsI}-\cite{CiarletLodsIV}, where it is shown that the solution of the linearized elasticity equations in a shell converges, when the shell thickness tends to zero, to different shell models, depending on the geometry of the shell and its boundary conditions. In particular, in the work of Ciarlet and Lods \cite{CiarletLodsIII}, the authors show that the Koiter’s shell model has the same asymptotic behavior, that is, its solutions converge, when the thickness of the shell tends to zero, to the same limit problems as the linearized elasticity equations do.

In this article we intend to justify a new two-dimensional flow model of a viscous fluid between two very close moving surfaces in a similar way to what was done in the above mentioned works \cite{CiarletLodsI}-\cite{CiarletLodsIII}, and, to confirm that when the distance between the two surfaces tends to zero, its behavior is the same as that observed in \cite{RodTabJMAA2021} for the Navier-Stokes equations, justifying it in sections \ref{sec-aa} and \ref{BoundaryConditions}, using a formal asymptotic expansion of the solution of the new model.

With this aim, in the first place, we will summarize, in section \ref{PreviousModels}, the results presented previously in our article \cite{RodTabJMAA2021}, which will allow us to make some assumptions about the behavior of the solutions of the Navier-Stokes. These hypothesis will be used in  sections \ref{NewAssumptions}-\ref{BoundaryConditions} to derive the new two-dimensional model proposed  in section \ref{NewModel}. Next, in section \ref{sec-aa}, we will begin the asymptotic analysis of the model, deriving in section \ref{subseccion-4-1} the limit model if the fluid velocity is known at the upper and lower bound surfaces, and obtaining, in section \ref{subseccion_sw}, the limit model when the tractions are known at the bound surfaces. Finally we will discuss the results achieved in section \ref{sec-conclusions}.

\section{Summary of the main previous results} \label{PreviousModels}

In our prior work \cite{RodTabJMAA2021} we studied the behavior of the Navier-Stokes equations in a domain bounded by two nearby moving surfaces, when the distance between them tends to zero. We observed that the asymptotic behavior of the solutions of the Navier-Stokes equations, in this case, strongly depends on the boundary conditions. In fact, two different limit models were obtained (one similar to a lubrication model and the other similar to a thin fluid layer model), depending on the boundary conditions in the original problem.

The two models presented in the preceding article \cite{RodTabJMAA2021} were derived from Navier-Stokes equations in  a three-dimensional thin domain, $\Omega^{\varepsilon}_t$, 
filled by a viscous fluid, that varies with time $t \in [0, T]$, given by
\begin{eqnarray}\Omega^{\varepsilon}_t&=&\left\{ (x_1^{\varepsilon},x_2^{\varepsilon},x_3^{\varepsilon})\in{R}^3:x_i(\xi_1,\xi_2,t)\leq x_i^{\varepsilon} \leq
x_i(\xi_1,\xi_2,t)+ h^\varepsilon(\xi_1,\xi_2,t)N_i(\xi_1,\xi_2,t),  \nonumber\right.\\
&&\left. (i=1,2,3), \ (\xi_1,\xi_2)\in D\subset \mathbb{R}^2
\right\} \label{eq-o-domain} \end{eqnarray} where
$\vec{X}_t(\xi_1,\xi_2)=\vec{X}(\xi_1,\xi_2,t)=(x_1(\xi_1,\xi_2,t),x_2(\xi_1,\xi_2,t),
x_3(\xi_1,\xi_2,t))$ is the lower
bound surface parametrization, $h^\varepsilon(\xi_1,\xi_2,t)$ is the gap between the two surfaces in
motion, and $\vec{N}(\xi_1,\xi_2,t)$ is the unit normal vector.

The lower bound surface is assumed to be regular
and the gap is assumed to be small with regard to the dimensions of the bound
surfaces. We take into account that the fluid film between the
surfaces is thin by introducing a small non-dimensional parameter 
$\varepsilon$, and setting that 
\begin{equation}
h^\varepsilon(\xi_1,\xi_2,t) = \varepsilon h(\xi_1,\xi_2,t), \quad
h(\xi_1,\xi_2,t) \ge h_0 > 0, \quad \forall \ (\xi_1,\xi_2)\in D\subset \mathbb{R}^2, \ \forall \ 
t\in [0,T].
\end{equation}

We introduce a reference domain 
\begin{equation}
\Omega=D \times [0,1] \label{eq-Omega}
\end{equation}
independent of
$\varepsilon$ and $t$, which is related to $\Omega^{\varepsilon}_t$ by the following change of variable: 
\begin{eqnarray}
    t^\varepsilon&=&t \label{eq-1-cv} \\
    x_i^\varepsilon &=& x_i(\xi_1,\xi_2,t)+\varepsilon \xi_3 h(\xi_1,\xi_2,t)N_i(\xi_1,\xi_2,t) \quad (i=1,2,3) \label{eq-2-cv}
\end{eqnarray}
where $(\xi_1,\xi_2)\in D$ and $\xi_3 \in[0,1]$. Now, given any scalar function $F^\varepsilon(t^\varepsilon,x_1^\varepsilon,x_2^\varepsilon,x_3^\varepsilon)$ defined on $\Omega^\varepsilon_t$, we can introduce another scalar function $F(\varepsilon)(t,\xi_1,\xi_2,\xi_3)$ on $\Omega$, using the change of variable:
\begin{equation}
F(\varepsilon)(t,\xi_1,\xi_2,\xi_3) = F^\varepsilon(t^\varepsilon,x_1^\varepsilon,x_2^\varepsilon,x_3^\varepsilon)
\end{equation}

We also define the basis $\left\{
\vec{a}_1,\vec{a}_2,\vec{a}_3\right\}$ 
\begin{eqnarray}
 \vec{a}_1(\xi_1,\xi_2,t)&=&\dfrac{\partial \vec{X}(\xi_1,\xi_2,t)}{\partial \xi_1} \label{base_a1} \\
 \vec{a}_2(\xi_1,\xi_2,t)&=&\dfrac{\partial \vec{X}(\xi_1,\xi_2,t)}{\partial \xi_2}\\
 \vec{a}_3(\xi_1,\xi_2,t)&=& \vec{N}(\xi_1,\xi_2,t) \label{base_a3}
\end{eqnarray}
so that, the velocity, $\vec{u}^{\varepsilon}$, and  the external density of 
volume forces, $\vec{f}^\varepsilon$, can be written in the new basis \eqref{base_a1}-\eqref{base_a3} as follows, where we adopt the convention of summing over repeated indices from 1 to 3, except where otherwise indicated:
\begin{eqnarray}\vec{u}^{\varepsilon} &=& u_i^{\varepsilon}
\vec{e}_i = u_k(\varepsilon)\vec{a}_{k}, \quad 
u_i^{\varepsilon}= \left ( u_k(\varepsilon)\vec{a}_{k} \right ) \cdot \vec{e}_i = u_k(\varepsilon) a_{ki} \label{cambio_base_u}\\
\vec{f}^\varepsilon &=& f_i^{\varepsilon}\vec{e}_i = 
f_k(\varepsilon)\vec{a}_{k}, \quad f_i^{\varepsilon}= \left ( f_k(\varepsilon)\vec{a}_{k} \right ) \cdot \vec{e}_i = f_k(\varepsilon){a}_{ki}
\label{cambio_base_f}
\end{eqnarray}
where ${a}_{ki} = \vec{a}_{k} \cdot \vec{e}_i$.

Taking into account \eqref{eq-1-cv}-\eqref{cambio_base_f}, Navier-Stokes equations can be written in the reference domain $\Omega$ in the following way (in the next equations, repeated indices indicate summation from 1 to 3, except for index $l$ and $n$, which take values from 1 to 2):
\begin{eqnarray}
&&\dfrac{ \partial u_k(\varepsilon)}{\partial t}  {a}_{ki} +
u_k(\varepsilon)\dfrac{
\partial {a}_{ki}}{\partial t} + \left({a}_{ki} \dfrac{
\partial u_k(\varepsilon)}{\partial \xi_n} + u_k(\varepsilon)\dfrac{ \partial {a}_{ki}}{\partial \xi_n} \right)\left[ -(\alpha_n \vec{a}_1 + \beta_n \vec{a}_2)\cdot\left(
\dfrac{\partial \vec{X}}{\partial t} + \varepsilon \xi_3 h
\dfrac{\partial \vec{a}_{3}}{\partial t} \right) \right]
\nonumber\\
&&\hspace*{+0.5cm}{}+\left({a}_{ki}\dfrac{
\partial u_k(\varepsilon)}{\partial \xi_3} + u_k(\varepsilon)\dfrac{
\partial {a}_{ki}}{\partial \xi_3} \right)\left(
-\dfrac{1}{\varepsilon h} \vec{a}_3 \cdot\dfrac{\partial
\vec{X}}{\partial t}- \dfrac{\xi_3}{ h} \dfrac{\partial h}{\partial
t}
\right)\nonumber\\
&&\hspace*{+0.5cm}{}+u_k(\varepsilon){a}_{kj} \left({a}_{qi} \dfrac{
\partial u_q(\varepsilon)}{\partial \xi_l} +u_q(\varepsilon)
\dfrac{
\partial {a}_{qi}}{\partial \xi_l} \right) \left(  \alpha_l
{a}_{1j} + \beta_l {a}_{2j}+ \gamma_l {a}_{3j}\right)\nonumber\\
&&\hspace*{+0.5cm}=-\dfrac{1}{\rho_0}\dfrac{
\partial p(\varepsilon)}{\partial \xi_l}\left(\alpha_l
{a}_{1i} + \beta_l {a}_{2i}+ \gamma_l {a}_{3i}\right) + \nu \left\{ \left[\dfrac{
\partial^2 (u_k(\varepsilon){a}_{ki})}{\partial \xi_l \partial \xi_m}
\left( \alpha_l {a}_{1j} + \beta_l {a}_{2j}+ \gamma_l {a}_{3j}
\right)\right.\right. \nonumber\\
&&\hspace*{+0.5cm}\left.\left.{}+ \dfrac{
\partial (u_k(\varepsilon){a}_{ki})}{\partial \xi_l}\dfrac{\partial}{\partial \xi_m}\left(
\alpha_l {a}_{1j} + \beta_l {a}_{2j}+ \gamma_l {a}_{3j} \right)
\right] \left( \alpha_m {a}_{1j} + \beta_m {a}_{2j}+ \gamma_m
{a}_{3j} \right) \right\} \nonumber\\
&&\hspace*{+0.5cm}{}+f_k(\varepsilon){a}_{ki}, \quad (i=1,2,3)\label{ec_ns_ij_alfa_beta}\\
&& \left({a}_{kj} \dfrac{
\partial u_k(\varepsilon)}{\partial \xi_l} + u_k(\varepsilon)\dfrac{
\partial {a}_{kj}}{\partial \xi_l}\right) \left(\alpha_l
{a}_{1j} + \beta_l {a}_{2j}+ \gamma_l {a}_{3j}\right) =0
\label{div_i_dr_alfa_beta}
\end{eqnarray}where $\alpha_l$, $\beta_l$ and $\gamma_l$ are defined in appendix \ref{ApendiceA} by expressions \eqref{alfaides}-\eqref{beta3n}.  We denote by $p$ the pressure, by $\rho_0$ the fluid density and by $\nu$ the kinematic viscosity. 

We begin assuming that $u_i(\varepsilon)$, $f_i(\varepsilon)$
($i=1,2,3$) and $p(\varepsilon)$ can be developed in powers of
$\varepsilon$, that is, 
\begin{eqnarray}
&& u_i(\varepsilon) = u_i^0 + \varepsilon u_i^1 + \varepsilon^2
u_i^2 + \cdots \quad (i=1,2,3) \label{ansatz_1}\\
&&p(\varepsilon) =\varepsilon^{-2} p^{-2} + \varepsilon^{-1} p^{-1}
+p^0 + \varepsilon p^1 + \varepsilon^2 p^2 + \cdots \label{ansatz_2}\\
&& f_i(\varepsilon) = f_i^0 + \varepsilon f_i^1 + \varepsilon^2
f_i^2 + \cdots \quad (i=1,2,3) \label{ansatz_3}
\end{eqnarray}

As mentioned above, using asymptotic analysis we are able to derive two different models depending on the boundary conditions chosen.

In the  first place, if we assume that the fluid slips at the lower surface $(\xi_3=0)$, and at the upper surface $(\xi_3=1)$, but there is continuity in the normal direction, so the tangential velocities at the lower and upper surfaces are known, and the normal velocity of each of them must match the fluid velocity, we obtain
\begin{eqnarray}
&&\hspace*{-0.9cm} \dfrac{1}{\sqrt{A^0}} \textrm{div}\left(\dfrac{h^3 }{ \sqrt{A^0}} M \nabla
p^{-2} \right) =12\mu \dfrac{\partial h}{\partial t} +
12\mu\dfrac{h A^1}{A^0} \left(\dfrac{\partial \vec{X}}{\partial t}
\cdot
\vec{a}_3\right)\nonumber\\
&&\hspace*{-0.5cm}{} - 6\mu \nabla h\cdot (\vec{W}^{0}-\vec{V}^{0}) +
\dfrac{6\mu h}{\sqrt{A^0}} \textrm{div} (\sqrt{A^0}(\vec{W}^{0} + \vec{V}^{0})) \label{Reynolds_gen}
\end{eqnarray}
that can be considered a generalization of Reynolds equation. We denote by $V_1 \vec{a}_1+V_2 \vec{a}_2$
the tangential velocity
at the lower surface and, by $W_1\vec{a}_1+W_2\vec{a}_2$ the
tangential velocity at the upper surface, and we have
\begin{eqnarray}
\vec{V}(\varepsilon)&=&(V_1,V_2)=\vec{V}^0+O(\varepsilon)\\
\vec{W}(\varepsilon)&=&(W_1,W_2)=\vec{W}^0+O(\varepsilon)
\end{eqnarray}   

Coefficients $A^0$, $A^1$ and matrix $M$ are defined in appendix \ref{ApendiceA} (\eqref{A0}-\eqref{M}), and $\mu = \rho_0 \nu$ is the dynamic viscosity.

Once obtained $p^{-2}$ using \eqref{Reynolds_gen}, the following approximation of the three components of the velocity is yielded
\begin{eqnarray}
&&\hspace*{-0.2cm} u_i^0 =\dfrac{h^2 (\xi_3^2-\xi_3)}{2\mu} \sum_{k=1}^2J^{0,0}_{ik} \dfrac{
\partial p^{-2}}{\partial \xi_k}+\xi_3(W_i^{0}-
V_i^{0})+ V_i^{0}, \quad (i=1,2)\label{ui_0_lub}\\
&&\hspace*{-0.2cm}
u_3^0  =  \dfrac{\partial \vec{X}}{\partial t} \cdot \vec{a}_3
\label{u30_lub}
\end{eqnarray}
where $J^{0,0}_{ik}$ is given by \eqref{J}.

If instead of considering that the tangential and normal velocities are known on the upper and lower surfaces, we assume that the normal component of the traction on $\xi_3=0$ and on $\xi_3=1$ are known pressures (denoted by $\pi_0^\varepsilon$ and $\pi_1^\varepsilon$, respectively), and that the tangential component of the traction on these surfaces are friction forces depending on the value of the velocities on $\partial D$, then we get a thin fluid layer model: 
\begin{eqnarray}
&& \hspace*{-0.5cm} u_i^0=W_i^{0}=V_i^{0} \quad (i=1,2) \label{u_i^0s}\\
&& \hspace*{-0.5cm} p^{-2}=p^{-1}=0\label{p-2p-1}\\
&& \hspace*{-0.5cm} p^0=\frac{2\mu}{h} \dfrac{\partial h}{\partial t} + \pi_0^0 \label{p0s-v2}
\\
&&\hspace*{-0.5cm}\dfrac{ \partial V_i^0}{\partial t}  + \sum_{l=1}^2 \left( V_l^0-C^0_l \right) \dfrac{
\partial V_i^0}{\partial \xi_l} +
\sum_{k=1}^2 
\left( R^0_{ik}+\sum_{l=1}^2
H^0_{ilk} V_l^0 \right) V_k^0=-\dfrac{1}{\rho_0} \sum_{l=1}^2 \dfrac{
\partial  \pi_0^0  }{\partial \xi_l}J^{0,0}_{il}  \nonumber \\ 
&&{} +\nu \left\{ \sum_{m=1}^2 \sum_{l=1}^2 
 \dfrac{
\partial^2 V_i^0 }{\partial \xi_m \partial \xi_l} J^{0,0}_{lm}
+ \sum_{k=1}^2 \sum_{l=1}^2  \dfrac{
\partial V_k^0 }{\partial \xi_l}\bar{L}^{0,0}_{ikl} 
+ \sum_{k=1}^2 V_k^0 \bar{S}_{ik}^{0,0}  + {\kappa}^0_i  \right\} \nonumber \\ 
&&{} + {F}^0_i-
 Q^0_{i3} \left ( \frac{\partial \vec{X}}{\partial t}
\cdot \vec{a}_3\right )  \quad (i=1,2)\label{ec_Vi0-v2} \end{eqnarray}
where coefficients $C^0_l$, $R^0_{ik}$, $H^0_{ilk}$, $J^{0,0}_{lm}$, $\bar{L}^{0,0}_{ikl}$, $\bar{S}_{ik}^{0,0}$,  ${\kappa}^0_i$, $F^0_i$ and $Q^0_{i3}$ are defined in appendix \ref{ApendiceA} (in \eqref{C}, \eqref{R}, \eqref{H}, \eqref{J}, \eqref{Lbar}, \eqref{Sbar},  \eqref{kappa}, \eqref{F} and \eqref{Q3} respectively),  and $\pi_0^0$ is the term of order zero on $\varepsilon$ of $\pi_0^\varepsilon$, that is,  $\pi_0^\varepsilon = \pi_0^0 + O(\varepsilon)$.

\begin{rmk}
Equation \eqref{ec_Vi0-v2} is exactly  the same as equation (168) in \cite{RodTabJMAA2021},  although some of the constants have been redefined here (see \eqref{Lbar} and \eqref{Sbar}) with respect to the definitions in \cite{RodTabJMAA2021} to simplify \eqref{ec_Vi0-v2}.
\end{rmk}

\section{New hypothesis about the dependence of the solution on $\xi_3$} \label{NewAssumptions}

If we carefully observe the steps of the proofs in the previous work \cite{RodTabJMAA2021}, we can see that $p^k\ (k=-2, -1, 0, 1)$, and $u_i^k\ (i=1,2,3; k=0,1)$ are polynomials in $\xi_3$ of at most degree three. Because of this, we are going to assume that, for $\varepsilon$  small enough, the following equalities are true:
\begin{align}
&\hspace*{-0.5cm} u_i(\varepsilon)(t,\xi_1,\xi_2,\xi_3)=\sum_{n=0}^3 \xi_3^n \bar{u}_i^n (\varepsilon) (t,\xi_1,\xi_2), \quad (i=1,2,3) \label{ui_pol_xi3}\\
&\hspace*{-0.5cm} p(\varepsilon)(t,\xi_1,\xi_2,\xi_3)=\sum_{n=0}^3 \xi_3^n \bar{p}^n(\varepsilon)(t,\xi_1,\xi_2),\label{p_pol_xi3} \\ &\hspace*{-0.5cm} f_i(\varepsilon)(t,\xi_1,\xi_2,\xi_3)=\sum_{n=0}^{\infty} \xi_3^n \bar{f}_i^n (\varepsilon) (t,\xi_1,\xi_2). \quad (i=1,2,3) \label{f_pol_xi3}\end{align}

We want to point out that the previous hypothesis is equivalent to neglecting in \eqref{ansatz_1}-\eqref{ansatz_2} the terms in $O(\varepsilon^2)$ when $\varepsilon$ is small.

Using expressions \eqref{ui_pol_xi3}-\eqref{f_pol_xi3} and \eqref{alfaides}-\eqref{gammades} (see appendix \ref{ApendiceA}), we can rewrite equations \eqref{ec_ns_ij_alfa_beta}-\eqref{div_i_dr_alfa_beta} as follows (repeated indices $k$, $j$ and $q$ indicate summation from 1 to 3, while repeated indices $l$ and $m$ indicate summation from 1 to 2):
\begin{eqnarray}
&&\hspace*{-0.5cm} \sum_{n=0}^3 \xi_3^n \dfrac{ \partial \bar{u}_k^n}{\partial t} \vec{a}_{k} +
\sum_{n=0}^3 \xi_3^n \bar{u}_k^n \dfrac{
\partial \vec{a}_{k}}{\partial t} \nonumber\\
&&{}- \left(\vec{a}_{k}\sum_{n=0}^3 \xi_3^n  \dfrac{
\partial \bar{u}_k^n}{\partial \xi_l} + \sum_{n=0}^3 \xi_3^n \bar{u}_k^n\dfrac{ \partial \vec{a}_{k}}{\partial \xi_l} \right)\left[ \sum_{r=0}^{\infty} (\varepsilon \xi_3 h)^r \left(\alpha_l^r \vec{a}_1 +  \beta_l^r \vec{a}_2\right) \cdot\left(
\dfrac{\partial \vec{X}}{\partial t} + \varepsilon \xi_3 h
\dfrac{\partial \vec{a}_{3}}{\partial t} \right) \right]
\nonumber\\
&&{}-\vec{a}_{k} \sum_{n=1}^3 n \xi_3^{n-1} \bar{u}_k^n\left[\dfrac{1}{\varepsilon h} \left(\vec{a}_3 \cdot\dfrac{\partial
\vec{X}}{\partial t} \right)+ \dfrac{\xi_3}{ h} \dfrac{\partial h}{\partial
t}\right.\nonumber\\
&&\left.{}+ \sum_{r=0}^{\infty} \varepsilon^{r} \xi_3^{r+1} h^{r-1} \left( \alpha_3^{r} \vec{a}_1 + \beta_3^{r} \vec{a}_2 \right) \cdot\left(
\dfrac{\partial \vec{X}}{\partial t} + \varepsilon \xi_3 h
\dfrac{\partial \vec{a}_{3}}{\partial t} \right)
\right]\nonumber\\
&&{}+\sum_{n=0}^3 \xi_3^n \bar{u}_k^n  \left(\vec{a}_{q}\sum_{d=0}^3 \xi_3^{d}  \dfrac{
\partial \bar{u}_q^{d}}{\partial \xi_l} + \sum_{d=0}^3 \xi_3^{d} \bar{u}_q^{d}\dfrac{ \partial \vec{a}_{q}}{\partial \xi_l} \right)   \sum_{r=0}^{\infty} (\varepsilon \xi_3 h)^{r} \left(\alpha_l^{r}
(\vec{a}_k \cdot \vec{a}_{1}) + \beta_l^{r} (\vec{a}_k \cdot \vec{a}_{2})\right))\nonumber\\
&&{}+\sum_{n=0}^3 \xi_3^n \bar{u}_k^n \left(\vec{a}_{q} \sum_{d=1}^3 d \xi_3^{d-1} \bar{u}_q^{d} \right) \left(  \sum_{r=0}^{\infty} \varepsilon^n \xi_3^{r+1} h^{r-1} \left( \alpha_3^{r}
(\vec{a}_k \cdot \vec{a}_{1}) + \beta_3^{r} (\vec{a}_k \cdot \vec{a}_{2})\right)\right.\nonumber\\
&&\left.{}+ \dfrac{1}{\varepsilon h} (\vec{a}_k \cdot \vec{a}_{3})\right)=-\dfrac{1}{\rho_0} \sum_{n=0}^3 \xi_3^n \dfrac{
\partial  \bar{p}^n}{\partial \xi_l} \sum_{r=0}^{\infty} (\varepsilon \xi_3 h)^{r} \left(\alpha_l^{r} \vec{a}_1 +  \beta_l^{r} \vec{a}_2\right) \nonumber\\
&&{} -\dfrac{1}{\rho_0} \sum_{n=1}^3 n \xi_3^{n-1} \bar{p}^n\left( \sum_{r=0}^{\infty} \varepsilon^{r} \xi_3^{{r}+1} h^{r-1} \left( \alpha_3^{r} \vec{a}_1 + \beta_3^{r} \vec{a}_2 \right) + \dfrac{1}{\varepsilon h} \vec{a}_{3}\right)\nonumber\\
&&{} + \nu \left[ \displaystyle\sum_{n=0}^3 \xi_3^n\dfrac{
\partial^2 (\vec{a}_{k}  \bar{u}_k^n)}{\partial \xi_l \partial \xi_m}
\sum_{r=0}^{\infty} (\varepsilon \xi_3 h)^{r} \left( \alpha_l^{r} {a}_{1j} + \beta_l^{r} {a}_{2j}
\right)\right. \nonumber\\
&&\left.{}+ 2\sum_{n=1}^3 n\xi_3^{n-1}\dfrac{
\partial (\vec{a}_{k}  \bar{u}_k^n)}{ \partial \xi_m}
\left(\sum_{r=0}^{\infty} \varepsilon^{r}
\xi_3^{r+1} h^{r-1} \left( \alpha_3^{r} {a}_{1j} + \beta_3^{r}
{a}_{2j}\right)+ \dfrac{1}{\varepsilon h}
{a}_{3j}
\right)\right. \nonumber\\
&&\left.{}+ \displaystyle\sum_{n=0}^3 \xi_3^n\dfrac{
\partial (\vec{a}_{k}  \bar{u}_k^n)}{\partial \xi_l}\dfrac{\partial}{\partial \xi_m} \left(\sum_{r=0}^{\infty} (\varepsilon \xi_3 h)^{r} \left( \alpha_l^{r} {a}_{1j} + \beta_l^{r} {a}_{2j}
\right) \right)\right. \nonumber\\
&&\left.{}+  \vec{a}_{k}
\displaystyle\sum_{n=1}^3 n \xi_3^{n-1}
\bar{u}_k^n\dfrac{\partial}{\partial \xi_m}\left(\sum_{r=0}^{\infty} \varepsilon^{r}
\xi_3^{r+1} h^{r-1} \left( \alpha_3^{r} {a}_{1j} + \beta_3^{r}
{a}_{2j}\right)+ \dfrac{1}{\varepsilon h}
{a}_{3j}  \right)
\right]\nonumber\\
&& \cdot \sum_{s=0}^{\infty} (\varepsilon \xi_3 h)^{s} \left( \alpha_m^{s} {a}_{1j} + \beta_m^{s} {a}_{2j}
\right) \nonumber\\
&&{}+ \nu \left[\vec{a}_{k} \sum_{n=2}^3 n(n-1) \xi_3^{n-2}
\bar{u}_k^n \left( \sum_{r=0}^{\infty} \varepsilon^{r}
\xi_3^{r+1} h^{r-1} \left( \alpha_3^{r} {a}_{1j} + \beta_3^{r}
{a}_{2j}\right)+ \dfrac{1}{\varepsilon h}
{a}_{3j}
\right)\right. \nonumber\\
&&\left.{}+ \displaystyle\sum_{n=0}^3 \xi_3^n \dfrac{
\partial (\vec{a}_{k} \bar{u}_k^n)}{\partial \xi_l} 
\sum_{r=1}^{\infty} r \varepsilon^r \xi_3^{r-1} h^{r} \left( \alpha_l^{r} {a}_{1j} + \beta_l^{r} {a}_{2j}
\right) \right. \nonumber\\
&&\left.{}+\vec{a}_{k} \displaystyle\sum_{n=1}^3 n \xi_3^{n-1}
\bar{u}_k^n \left( \sum_{r=0}^{\infty} (r+1)\varepsilon^{r} \xi_3^{r}
h^{r-1} \left( \alpha_3^{r} {a}_{1j} + \beta_3^{r} {a}_{2j}
\right) \right)\right]\nonumber\\
&& \cdot\left( \sum_{s=0}^{\infty} \varepsilon^{s}
\xi_3^{s+1} h^{s-1} \left( \alpha_3^{s} {a}_{1j} + \beta_3^{s}
{a}_{2j}\right) + \dfrac{1}{\varepsilon h}
{a}_{3j} \right)+
\sum_{n=0}^{\infty} \xi_3^n \bar{f}_k^n \vec{a}_{k}\label{ec_ns_ij_alfa_beta_pol_xi3}
\\
&&\hspace*{-0.5cm}\sum_{n=0}^3 \xi_3^n \dfrac{
\partial \bar{u}_k^n}{\partial \xi_l} \sum_{r=0}^{\infty} (\varepsilon \xi_3 h)^{r} \left( \alpha_l^{r}
(\vec{a}_k \cdot \vec{a}_{1}) + \beta_l^{r} (\vec{a}_k \cdot \vec{a}_{2})\right) \nonumber\\
&&{} +\sum_{n=0}^3 \xi_3^n \bar{u}_k^n\dfrac{
\partial {a}_{kj}}{\partial \xi_l} \sum_{r=0}^{\infty} (\varepsilon \xi_3 h)^{r} \left( \alpha_l^{r} {a}_{1j} + \beta_l^{r} {a}_{2j}
\right) \nonumber\\
&&{}+ \sum_{n=1}^3 n \xi_3^{n-1 } \bar{u}_k^n \left(  \sum_{r=0}^{\infty} \varepsilon^{r}
\xi_3^{r+1} h^{r-1} \left( \alpha_3^{r}
(\vec{a}_k \cdot \vec{a}_{1}) + \beta_3^{r} (\vec{a}_k \cdot \vec{a}_{2})\right) + \dfrac{1}{\varepsilon h} (\vec{a}_k \cdot \vec{a}_{3})\right) =0
\label{div_i_dr_alfa_beta_pol_xi3}
\end{eqnarray}
and identify the terms multiplied by  $\xi_3^n$ ($n=0,1,2,3$) in  \eqref{ec_ns_ij_alfa_beta_pol_xi3}-\eqref{div_i_dr_alfa_beta_pol_xi3}.  In equation \eqref{ec_ns_xi3^n}, below,  repeated indices $k$ and $q$ indicate, again, summation from 1 to 3, while repeated indices $l$ and $m$ indicate summation from 1 to 2. 
\begin{eqnarray}
&&\hspace*{-0.5cm}  \dfrac{ \partial \bar{u}_k^n}{\partial t} \vec{a}_{k} +
 \bar{u}_k^n \dfrac{
\partial \vec{a}_{k}}{\partial t} -  \left(\vec{a}_{k} \dfrac{
\partial \bar{u}_k^n}{\partial \xi_l} +  \bar{u}_k^n\dfrac{ \partial \vec{a}_{k}}{\partial \xi_l} \right) C^0_l\nonumber\\
&&{}- \sum_{m=0}^{n-1} \left(\vec{a}_{k} \dfrac{
\partial \bar{u}_k^m}{\partial \xi_l} +  \bar{u}_k^m\dfrac{ \partial \vec{a}_{k}}{\partial \xi_l} \right)  (\varepsilon  h)^{n-m} C^{n-m,n-m-1}_l \nonumber\\
&&{}-\dfrac{n+1}{\varepsilon h} \vec{a}_{k} \bar{u}_k^{n+1}\left(\vec{a}_3 \cdot\dfrac{\partial
\vec{X}}{\partial t} \right)- \dfrac{n}{ h}  \vec{a}_{k}   \bar{u}_k^n \left(\dfrac{\partial h}{\partial
t} + C^0_3\right)\nonumber\\
&&{} -\vec{a}_{k} \sum_{m=0}^{n-2} (m+1) \bar{u}_k^{m+1} \varepsilon^{n-m-1} h^{n-m-2} C_3^{n-m-1,n-m-2}
\nonumber\\
&&{}+\sum_{m=0}^n \bar{u}_k^m \sum_{j=0}^{n-m}  \left(\vec{a}_{q} \dfrac{
\partial \bar{u}_q^j}{\partial \xi_l} + \bar{u}_q^j\dfrac{ \partial \vec{a}_{q}}{\partial \xi_l} \right)  (\varepsilon  h)^{n-m-j} B_{lk}^{n-m-j}
\nonumber\\
&&{}+ \sum_{m=0}^{n-1}  \bar{u}_k^m  \sum_{j=1}^{n-m} j \varepsilon^{n-m-j} h^{n-m-j-1} B_{3k}^{n-m-j} \vec{a}_{q}
 \bar{u}_q^{j} +\dfrac{1}{\varepsilon h} \sum_{m=0}^n \bar{u}_3^m (n-m+1) \left(\vec{a}_{q}   \bar{u}_q^{n-m+1} \right)\nonumber\\
&&{}   =-\dfrac{1}{\rho_0} \sum_{m=0}^n \dfrac{
\partial  \bar{p}^m}{\partial \xi_l}  (\varepsilon  h)^{n-m} \left(\alpha_l^{n-m} \vec{a}_1 +  \beta_l^{n-m} \vec{a}_2\right) - \dfrac{n+1}{\varepsilon h \rho_0} \bar{p}^{n+1} \vec{a}_{3}\nonumber\\
&&{} -\dfrac{1}{\rho_0} \sum_{m=1}^{n}   m  \bar{p}^{m}  \varepsilon^{n-m}  h^{n-m-1} \left( \alpha_3^{n-m} \vec{a}_1 + \beta_3^{n-m} \vec{a}_2 \right) \nonumber\\
&&{} + \nu \left[ \displaystyle\sum_{r=0}^n (\varepsilon  h)^{n-r} \dfrac{
\partial^2 (\vec{a}_{k}  \bar{u}_k^r)}{\partial \xi_l \partial \xi_m}
\sum_{s=0}^{n-r} J_{l,m}^{s,n-r-s}  + 2 \sum_{r=1}^n r \varepsilon^{n-r}  h^{n-r-1} \dfrac{
\partial (\vec{a}_{k}  \bar{u}_k^r)}{ \partial \xi_m}
\sum_{s=0}^{n-r} 
J_{3m}^{s,n-r-s}
\right. \nonumber\\
&&\left.{}+ \displaystyle\sum_{r=0}^n (\varepsilon  h)^{n-r}\dfrac{
\partial (\vec{a}_{k}  \bar{u}_k^r)}{\partial \xi_l} \sum_{s=0}^{n-r} K_l^{s,n-r-s} + \displaystyle\sum_{r=0}^n \varepsilon^{n-r}  h^{n-r-1} \dfrac{
\partial (\vec{a}_{k}  \bar{u}_k^r)}{\partial \xi_l} \sum_{s=1}^{n-r} s \dfrac{\partial h}{\partial \xi_m} J_{lm}^{s,n-r-s}\right. \nonumber\\
&&\left.{}+  \vec{a}_{k}
\displaystyle\sum_{r=1}^n r  \varepsilon^{n-r}  h^{n-r-1} 
\bar{u}_k^r \sum_{s=0}^{n-r} K_3^{s,n-r-s}+  \vec{a}_{k}
\displaystyle\sum_{r=1}^n r \varepsilon^{n-r} h^{n-r-2}
\bar{u}_k^r \sum_{s=0}^{n-r} (s-1) \dfrac{\partial h}{\partial \xi_m}  J_{3m}^{s,n-r-s} 
  \right. \nonumber\\
&&\left.{}+  \vec{a}_{k}
\sum_{r=0}^n (r+1) \varepsilon^{n-r-1} h^{n-r-1} 
\bar{u}_k^{r+1} H^{n-r}_{mm3}
\right]\nonumber\\
&&{}+ \nu \left[\vec{a}_{k} \sum_{r=2}^{n} r(r-1)  \varepsilon^{n-r}  h^{n-r-2}
\bar{u}_k^r \sum_{s=0}^{n-r} 
J_{33}^{s,n-r-s} 
+ \dfrac{\vec{a}_{k}}{\varepsilon^2 h^2} (n+2)(n+1) \bar{u}_k^{n+2} \right. \nonumber\\
&&\left.{}+ \displaystyle\sum_{r=0}^{n-1}  \varepsilon^{n-r}  h^{n-r-1}\dfrac{
\partial (\vec{a}_{k} \bar{u}_k^r)}{\partial \xi_l} 
\sum_{s=1}^{n-r} s   J_{l3}^{s,n-r-s} 
+\vec{a}_{k} \displaystyle\sum_{r=1}^n r \varepsilon^{n-r}  h^{n-r-2}  \bar{u}_k^r
 \sum_{s=0}^{n-r} (s+1) J_{33}^{s,n-r-s} \right]\nonumber\\
&& {} +
 \bar{f}_k^n \vec{a}_{k}, \quad (n=0,1,2,3) \label{ec_ns_xi3^n}
\\
&&\hspace*{-0.5cm}\bar{u}_3^{n+1}  =- \dfrac{\varepsilon h}{n+1} \sum_{m=0}^n  (\varepsilon  h)^{n-m} \sum_{k=1}^2 \left[ \sum_{l=1}^2  \left(\dfrac{
\partial \bar{u}_k^m}{\partial \xi_l}   B_{lk}^{n-m}  +  \bar{u}_k^m H_{llk}^{n-m} \right)+   \bar{u}_3^m    H_{kk3}^{n-m}
 +  \dfrac{ m}{h}  \bar{u}_k^m    
B_{3k}^{n-m}  \right]\nonumber\\
&& (n=0,1,2)
\label{div_u3n}\\
&&\hspace*{-0.5cm} \sum_{m=0}^3  (\varepsilon h)^m \left[
  \sum_{k=1}^2  \sum_{l=1}^2 \dfrac{
\partial \bar{u}_k^{3-m}}{\partial \xi_l}   B^m_{lk} +  \sum_{k=1}^3 \bar{u}_k^{3-m} \sum_{l=1}^2 H^m_{llk}  + \dfrac{3-m}{h}  \sum_{k=1}^2 \bar{u}_k^{3-m} B^m_{3k}\right]=0
\label{ec_div_n3}
\end{eqnarray} 
where we have introduced the notation 
$\bar{u}_i^4=\bar{u}_i^5=0$, \ ($i=1,2,3$). The coefficients $B^j_{lk}$, $C^0_l$, $C^{i,j}_l$, $H^j_{ilk}$, $J^{i,j}_{lm}$, $K^{j,i}_l$ are given by \eqref{B}, \eqref{C}, \eqref{Cij}, \eqref{H}, \eqref{J}, \eqref{Kji}.

We multiply equations \eqref{ec_ns_xi3^n}  by $\alpha^0_i \vec{a}_i$ ($i=1,2$) and sum in $i$, then we repeat the procedure multiplying \eqref{ec_ns_xi3^n} by $\beta^0_i \vec{a}_i$ ($i=1,2$) and adding in $i$ again, to yield these  equations: 
\begin{eqnarray}
&&\hspace*{-0.5cm}  \dfrac{ \partial \bar{u}_i^n}{\partial t}  +
\sum_{k=1}^3 \bar{u}_k^n Q^0_{ik} -  \sum_{l=1}^2 \dfrac{
\partial \bar{u}_i^n}{\partial \xi_l}   C^0_l  - \sum_{m=0}^{n-1} \sum_{l=1}^2  \left(\dfrac{
\partial \bar{u}_i^m}{\partial \xi_l} + \sum_{k=1}^3 \bar{u}_k^m H^0_{ilk} \right)  (\varepsilon  h)^{n-m} C^{n-m,n-m-1}_l \nonumber\\
&&{}-\dfrac{n+1}{\varepsilon h} \bar{u}_i^{n+1}\left(\vec{a}_3 \cdot\dfrac{\partial
\vec{X}}{\partial t} \right)- \dfrac{n}{ h}  \bar{u}_i^n \left(\dfrac{\partial h}{\partial
t} + C^0_3\right)\nonumber\\
&&{} -\sum_{m=0}^{n-2} (m+1) \bar{u}_i^{m+1} \varepsilon^{n-m-1} h^{n-m-2} C_3^{n-m-1,n-m-2}
\nonumber\\
&&{}+\sum_{m=0}^n \sum_{k=1}^2 \bar{u}_k^m \sum_{j=0}^{n-m} \sum_{l=1}^2 \left( \dfrac{
\partial \bar{u}_i^j}{\partial \xi_l} + \sum_{q=1}^3 \bar{u}_q^j H^0_{ilq} \right)  (\varepsilon  h)^{n-m-j} B_{lk}^{n-m-j}
\nonumber\\
&&{}+ \sum_{m=0}^{n-1} \sum_{k=1}^2 \bar{u}_k^m  \sum_{j=1}^{n-m} j \varepsilon^{n-m-j} h^{n-m-j-1} B_{3k}^{n-m-j}  \bar{u}_i^{j} +\dfrac{1}{\varepsilon h} \sum_{m=0}^n \bar{u}_3^m (n-m+1) \bar{u}_i^{n-m+1} \nonumber\\
&&{}   =-\dfrac{1}{\rho_0} \sum_{m=0}^n \dfrac{
\partial  \bar{p}^m}{\partial \xi_l}  (\varepsilon  h)^{n-m} J_{il}^{0,n-m}  -\dfrac{1}{\rho_0} \sum_{m=1}^{n}   m  \bar{p}^{m}  \varepsilon^{n-m}  h^{n-m-1} J_{i3}^{0,n-m} \nonumber\\
&&{} + \nu \left\{ \sum_{r=0}^n \varepsilon^{n-r} \sum_{m=1}^2\sum_{l=1}^2 \dfrac{
\partial^2   \bar{u}_i^r}{\partial \xi_l \partial \xi_m}
\iota_{lm}^{n-r}+ \sum_{r=0}^n \varepsilon^{n-r}   \sum_{k=1}^3 \sum_{l=1}^2 \dfrac{
\partial  \bar{u}_k^r}{\partial \xi_l} L_{ikl}^{n,r} + \sum_{r=0}^n \varepsilon^{n-r} \sum_{k=1}^3 \bar{u}_k^r S_{ik}^{n,r}
\right\}\nonumber\\
&&{}+  
\dfrac{\nu(n+1)}{\varepsilon h} 
\bar{u}_i^{n+1} \sum_{m=1}^2 H^{0}_{mm3}+ \nu  
\dfrac{ (n+2)(n+1) }{\varepsilon^2 h^2}\bar{u}_i^{n+2} +
 \bar{f}_i^n, \quad (i=1,2, \ n=0,1,2,3) \label{ec_ns_xi3^n_ai}\end{eqnarray} where $Q^0_{ik}$, $\iota^n_{lm}$, $L^{n,r}_{ikl}$ and $S^{n,r}_{ik}$ are given by \eqref{Q}, \eqref{Q3}, \eqref{iotalm}, \eqref{Likl_nr} and \eqref{Sik_nr} respectively.

If equations \eqref{ec_ns_xi3^n} are multiplied by $\vec{a}_3$, we obtain: 
\begin{eqnarray}
&&\hspace*{-0.5cm}  \dfrac{ \partial \bar{u}_3^n}{\partial t}  +\sum_{k=1}^2
 \bar{u}_k^n \left[\dfrac{
\partial \vec{a}_{k}}{\partial t} \cdot \vec{a}_3 -\sum_{l=1}^2 C^0_l \dfrac{ \partial \vec{a}_{k}}{\partial \xi_l} \cdot \vec{a}_3\right]- \sum_{l=1}^2 \dfrac{
\partial \bar{u}_3^n}{\partial \xi_l} C^0_l\nonumber\\
&&{}- \sum_{m=0}^{n-1}\sum_{l=1}^2  \left[ \dfrac{
\partial \bar{u}_3^m}{\partial \xi_l} + \sum_{k=1}^{2} \bar{u}_k^m \left(\dfrac{ \partial \vec{a}_{k}}{\partial \xi_l} \cdot \vec{a}_3 \right)\right]  (\varepsilon  h)^{n-m} C^{n-m,n-m-1}_l \nonumber\\
&&{}-\dfrac{n+1}{\varepsilon h}  \bar{u}_3^{n+1}\left(\vec{a}_3 \cdot\dfrac{\partial
\vec{X}}{\partial t} \right)- \dfrac{n}{ h}    \bar{u}_3^n \left(\dfrac{\partial h}{\partial
t} + C^0_3\right)\nonumber\\
&&{} - \sum_{m=0}^{n-2} (m+1) \bar{u}_3^{m+1} \varepsilon^{n-m-1} h^{n-m-2} C_3^{n-m-1,n-m-2}
\nonumber\\
&&{}+\sum_{m=0}^n \sum_{k=1}^2\bar{u}_k^m \sum_{j=0}^{n-m} \sum_{l=1}^2 \left( \dfrac{
\partial \bar{u}_3^j}{\partial \xi_l} + \bar{u}_q^j\dfrac{ \partial \vec{a}_{q}}{\partial \xi_l} \cdot \vec{a}_3 \right)  (\varepsilon  h)^{n-m-j} B_{lk}^{n-m-j}
\nonumber\\
&&{}+ \sum_{m=0}^{n-1} \sum_{k=1}^2 \bar{u}_k^m  \sum_{j=0}^{n-m-1} (n-m-j)\varepsilon^{j} h^{j-1} B_{3k}^j 
 \bar{u}_3^{n-m-j} +\dfrac{1}{\varepsilon h} \sum_{m=0}^n \bar{u}_3^m (n-m+1)  \bar{u}_3^{n-m+1} \nonumber\\
&&{}   = - \dfrac{n+1}{\varepsilon h \rho_0} \bar{p}^{n+1} + \nu \left[ \sum_{r=0}^n \varepsilon^{n-r} \sum_{l=1}^2\sum_{m=1}^2\dfrac{
\partial^2  \bar{u}_3^r}{\partial \xi_l \partial \xi_m} 
\iota_{lm}^{n-r}+\sum_{r=0}^n \varepsilon^{n-r} \sum_{k=1}^3 \sum_{l=1}^2 \dfrac{
\partial  \bar{u}_k^r}{\partial \xi_l} L_{3kl}^{n,r} \right. \nonumber\\
&&\left.{}+\sum_{r=0}^n \varepsilon^{n-r}  \sum_{k=1}^2\bar{u}_k^r S_{3k}^{n,r}
\right]+  
 \dfrac{\nu (n+1)}{\varepsilon h}
\bar{u}_3^{n+1} \sum_{m=1}^2 H^{0}_{mm3}+ \dfrac{\nu}{\varepsilon^2 h^2} (n+2)(n+1) \bar{u}_3^{n+2} \nonumber\\
&&{}+
 \bar{f}_3^n, \quad (n=0,1,2,3) \label{ec_ns_xi3^n_a3}
\end{eqnarray} 
where the coefficients $L^{n,r}_{3kl}$ and $S^{n,r}_{3k}$ are defined in \eqref{L3kl_nr} and \eqref{S3k_nr}.

Since we have assumed that the velocity and the pressure are polynomials of degree three in $\xi_3$ (\eqref{ui_pol_xi3}-\eqref{p_pol_xi3}), we have 16 unknowns to determine. Out of these unknowns,  the terms $\bar{u}^k_3$ and $\bar{p}^k$ ($k=1,2,3$) corresponding to the third component of the velocity and the pressure, respectively, are given by \eqref{div_u3n} and \eqref{ec_ns_xi3^n_a3} using the terms $\bar{u}^k_i$ ($i=1,2, k=0,1,2$) once they have been computed. Therefore, we must actually determine 10 unknowns.

Denoting, as previously done, by  $V_1 \vec{a}_1+V_2 \vec{a}_2$ and $W_1\vec{a}_1+W_2\vec{a}_2$ the tangential velocity
at the lower and upper surfaces, respectively, we have
\begin{eqnarray}
u_k^{\varepsilon} \vec{e}_k = u_k(\varepsilon)\vec{a}_{k}&=&
V_1(\varepsilon) \vec{a}_1+V_2(\varepsilon) \vec{a}_2+\left( \dfrac{\partial
\vec{X}}{\partial t}
\cdot \vec{a}_3\right)\vec{a}_3 \ \textrm{on }\xi_3=0 \label{cc_xi3_0}\\
u_k^{\varepsilon} \vec{e}_k =u_k(\varepsilon)\vec{a}_{k}&=&
W_1(\varepsilon) \vec{a}_1+ W_2(\varepsilon) \vec{a}_2+\left( \dfrac{\partial (\vec{X}+
\varepsilon h\vec{a}_3)}{\partial t} \cdot \vec{a}_3\right)\vec{a}_3
\ \textrm{on }\xi_3=1 \label{cc_xi3_1}
\end{eqnarray}
and, taking into account \eqref{ui_pol_xi3}, we yield
\begin{eqnarray}
&&\bar{u}_i^0=V_i \quad (i=1,2) \label{ui0_Vi}\\
&&\bar{u}_3^0= \dfrac{\partial \vec{X}}{\partial t}
\cdot \vec{a}_3 \label{u30}\\
&&\sum_{k=1}^3  \bar{u}_i^k= W_i - V_i \quad (i=1,2)
 \label{sum_uik_Wi_Vi}\\
&&\sum_{k=1}^3  \bar{u}_3^k=  \varepsilon \dfrac{\partial
h}{\partial t}
\label{sum_u3k_h}
\end{eqnarray} 

Equality \eqref{u30} gives us an expression for $\bar{u}_3^0$, so it is no longer an unknown, it is determined by the lower bound surface. At this point, 9 unknowns are left, $\bar{u}_i^k$ $(i=1,2,\ k=0,1,2,3)$ and $\bar{p}^0$, but we will see that not all are needed to obtain an approximation of the velocity and the pressure.

\section{New Model} \label{NewModel}

As we have just seen in section \ref{NewAssumptions}, operating with the truncated part of the formal asymptotic expansion of the solution,  hypotheses \eqref{ui_pol_xi3}-\eqref{f_pol_xi3} allow us to derive a two-dimensional limit model, which we shall call new model from now on, formed by equations:
\begin{align}
&\hspace*{-0.5cm} u_i(\varepsilon)(t,\xi_1,\xi_2,\xi_3)=\sum_{n=0}^3 \xi_3^n \bar{u}_i^n (\varepsilon) (t,\xi_1,\xi_2), \quad (i=1,2,3) \label{ui_pol_xi3m}\\
&\hspace*{-0.5cm} p(\varepsilon)(t,\xi_1,\xi_2,\xi_3)=\sum_{n=0}^3 \xi_3^n \bar{p}^n(\varepsilon)(t,\xi_1,\xi_2),\label{p_pol_xi3m} \\ 
&\hspace*{-0.5cm}  \dfrac{ \partial \bar{u}_i^n}{\partial t}  +
\sum_{k=1}^3 \bar{u}_k^n Q^0_{ik} -  \sum_{l=1}^2 \dfrac{
\partial \bar{u}_i^n}{\partial \xi_l}   C^0_l - \dfrac{n}{ h}  \bar{u}_i^n \left(\dfrac{\partial h}{\partial
t} + C^0_3\right) \nonumber\\
&{}+\sum_{m=0}^n \sum_{k=1}^2 \bar{u}_k^m \sum_{j=0}^{n-m} \sum_{l=1}^2 \left( \dfrac{
\partial \bar{u}_i^j}{\partial \xi_l} + \sum_{q=1}^3 \bar{u}_q^j H^0_{ilq} \right)  (\varepsilon  h)^{n-m-j} B_{lk}^{n-m-j}
\nonumber\\
&{}+ \sum_{m=0}^{n-1} \sum_{k=1}^2 \bar{u}_k^m  \sum_{j=1}^{n-m} j \varepsilon^{n-m-j} h^{n-m-j-1} B_{3k}^{n-m-j} 
 \bar{u}_i^{j} +\dfrac{1}{\varepsilon h} \sum_{m=1}^n \bar{u}_3^m (n-m+1) \bar{u}_i^{n-m+1} \nonumber\\
&{} - \sum_{m=0}^{n-1} \sum_{l=1}^2  \left(\dfrac{
\partial \bar{u}_i^m}{\partial \xi_l} + \sum_{k=1}^3 \bar{u}_k^m H^0_{ilk} \right)  (\varepsilon  h)^{n-m} C^{n-m,n-m-1}_l\nonumber\\
&{} -\sum_{m=0}^{n-2} (m+1) \bar{u}_i^{m+1} \varepsilon^{n-m-1} h^{n-m-2} C_3^{n-m-1,n-m-2}
\nonumber\\
&{}   =-\dfrac{1}{\rho_0} \sum_{m=0}^n \sum_{l=1}^2\dfrac{
\partial  \bar{p}^m}{\partial \xi_l}  (\varepsilon  h)^{n-m} J_{il}^{0,n-m}  -\dfrac{1}{\rho_0} \sum_{m=1}^{n}   m  \bar{p}^{m}  \varepsilon^{n-m}  h^{n-m-1} J_{i3}^{0,n-m} \nonumber\\
&{} + \nu \sum_{r=0}^n \varepsilon^{n-r}  \left[ \sum_{m=1}^2\sum_{l=1}^2 \dfrac{
\partial^2   \bar{u}_i^r}{\partial \xi_l \partial \xi_m}
\iota_{lm}^{n-r}+    \sum_{k=1}^3 \sum_{l=1}^2 \dfrac{
\partial  \bar{u}_k^r}{\partial \xi_l} L_{ikl}^{n,r} +  \sum_{k=1}^3 \bar{u}_k^r S_{ik}^{n,r}
\right] \nonumber\\
&{}+  
\dfrac{\nu(n+1)}{\varepsilon h} \dfrac{A^1}{A^0}
\bar{u}_i^{n+1} + \nu  
 {\dfrac{ (n+2)(n+1) }{\varepsilon^2 h^2}\bar{u}_i^{n+2}} +
 \bar{f}_i^n, \quad (i=1,2, n=0,1,2,3) \label{ec_uin}
\\
&\hspace*{-0.5cm} \sum_{m=0}^3  (\varepsilon h)^m \left[
  \sum_{k=1}^2  \sum_{l=1}^2 \dfrac{
\partial \bar{u}_k^{3-m}}{\partial \xi_l}   B^m_{lk} +  \sum_{k=1}^3 \bar{u}_k^{3-m} \sum_{l=1}^2 H^m_{llk}  + \dfrac{3-m}{h}  \sum_{k=1}^2 \bar{u}_k^{3-m} B^m_{3k}\right]=0
\label{ec_div_n3m}\\
&\hspace*{-0.5cm}\bar{u}_3^0= \dfrac{\partial \vec{X}}{\partial t}
\cdot \vec{a}_3 \label{u30m}\\
&\hspace*{-0.5cm}\bar{u}_3^{n+1}   =- \dfrac{\varepsilon h}{n+1} \sum_{m=0}^n  (\varepsilon  h)^{n-m} \sum_{k=1}^2 \left[ \sum_{l=1}^2  \left(\dfrac{
\partial \bar{u}_k^m}{\partial \xi_l}   B_{lk}^{n-m}  +  \bar{u}_k^m H_{llk}^{n-m} \right)+   \bar{u}_3^m    H_{kk3}^{n-m}
 +  \dfrac{ m}{h}  \bar{u}_k^m    
B_{3k}^{n-m}  \right] \nonumber\\
&(n=0,1,2)
\label{u3n}\\
&\hspace*{-0.5cm}  \bar{p}^{n+1}= {\dfrac{\mu}{\varepsilon h} (n+2)\bar{u}_3^{n+2}} +  
\dfrac{\mu A^1}{A^0} 
\bar{u}_3^{n+1} - \dfrac{ \rho_0}{n+1} \sum_{m=1}^n \bar{u}_3^m (n-m+1)  \bar{u}_3^{n-m+1} \nonumber\\
&{} +\dfrac{\varepsilon h \rho_0}{n+1} \left\{- \dfrac{ \partial \bar{u}_3^n}{\partial t}  -\sum_{k=1}^2
 \bar{u}_k^n Q^0_{3k}+ \sum_{l=1}^2 \dfrac{
\partial \bar{u}_3^n}{\partial \xi_l} C^0_l + \dfrac{n}{ h}    \bar{u}_3^n \left(\dfrac{\partial h}{\partial
t} + C^0_3\right) \right.\nonumber\\
&{} + \sum_{m=0}^{n-1} \sum_{l=1}^2 \left[ \dfrac{
\partial \bar{u}_3^m}{\partial \xi_l} + \sum_{k=1}^2 \bar{u}_k^m \left(\dfrac{ \partial \vec{a}_{k}}{\partial \xi_l} \cdot \vec{a}_3\right) \right]  (\varepsilon  h)^{n-m} C^{n-m,n-m-1}_l \nonumber\\
&{} +\sum_{m=0}^{n-2} (m+1) \bar{u}_3^{m+1} \varepsilon^{n-m-1} h^{n-m-2} C_3^{n-m-1,n-m-2}
\nonumber\\
&{}-\sum_{m=0}^n \sum_{k=1}^2\bar{u}_k^m \sum_{j=0}^{n-m} \sum_{l=1}^2 \left[ \dfrac{
\partial \bar{u}_3^j}{\partial \xi_l} + \bar{u}_q^j \left(\dfrac{ \partial \vec{a}_{q}}{\partial \xi_l} \cdot \vec{a}_3 \right)\right]  (\varepsilon  h)^{n-m-j} B_{lk}^{n-m-j}
\nonumber\\
&{}- \sum_{m=0}^{n-1} \sum_{k=1}^2 \bar{u}_k^m  \sum_{j=0}^{n-m-1} (n-m-j)\varepsilon^{j} h^{j-1} B_{3k}^j 
 \bar{u}_3^{n-m-j}\nonumber\\
&\left.{}   + \nu \sum_{r=0}^n \varepsilon^{n-r}  \left[ \sum_{l=1}^2\sum_{m=1}^2\dfrac{
\partial^2  \bar{u}_3^r}{\partial \xi_l \partial \xi_m} 
\iota_{lm}^{n-r}+ \sum_{k=1}^3 \sum_{l=1}^2 \dfrac{
\partial  \bar{u}_k^r}{\partial \xi_l} L_{3kl}^{n,r} +  \sum_{k=1}^3\bar{u}_k^r S_{3k}^{n,r}
\right]+
 \bar{f}_3^n\right\} \nonumber\\& (n=0,1,2) \label{pn}
\end{align}

Examining the new model, we observe that the equations can be divided into  two groups: a first group, including equations \eqref{ec_uin}  and \eqref{ec_div_n3m}, that must be solved to obtain the terms  $\bar{u}_1^0$, $\bar{u}_2^0$, $\bar{p}^0$, $\bar{u}_1^1$, $\bar{u}_2^1$, $\bar{u}_1^2$, $\bar{u}_2^2$, $\bar{u}_1^3$ and $\bar{u}_2^3$, and a second group, including equations \eqref{u30m}-\eqref{pn}, that allow us to eliminate the terms $\bar{u}_3^0$, $\bar{u}_3^1$, $\bar{u}_3^2$, $\bar{u}_3^3$, $\bar{p}^1$, $\bar{p}^2$, $\bar{p}^3$ from equations \eqref{ec_uin}-\eqref{ec_div_n3m}. Once the aforementioned elimination has been carried out, we can solve the first group of equations to compute $\bar{u}_i^k$ $(k=0,1,2,3, i=1,2)$ and $\bar{p}^0$, and use the second group of equations to obtain $\bar{u}_3^k$ ($k=0,1,2,3$) and $\bar{p}^j$ ($j=1,2,3$). Boundary and initial conditions must be added to this system of equations. 

\section{Asymptotic Analysis of the new model} \label{sec-aa} 

This section is devoted to the asymptotic analysis of the model \eqref{ui_pol_xi3m}-\eqref{pn}, proposed in the previous section. We want to check if the asymptotic behavior of the new model, when $\varepsilon$ tends to zero, is the same as the Navier-Stokes equations, shown in \cite{RodTabJMAA2021}.

Let us start assuming that $\bar{u}_i^n(\varepsilon)$, $\bar{f}_i^n(\varepsilon)$, $\bar{p}^n(\varepsilon)$, $W_i(\varepsilon)$ and $V_i(\varepsilon)$ ($i=1,2,3$, $n=0,1,2,3$) can be developed in powers of
$\varepsilon$, that is:
\begin{eqnarray}
&& \bar{u}_i^n(\varepsilon) = \sum_{k=0}^{\infty} \varepsilon^k \bar{u}_i^{n,k}   \quad (i=1,2,3, \ n=0,1,2,3) \label{ansatz_bar_u}\\
&&\bar{p}^n(\varepsilon) = \sum_{k=-2}^{\infty}\varepsilon^{k} \bar{p}^{n,k} \quad (n=0,1,2,3) \label{ansatz_bar_p}\\
&& \bar{f}_i^n(\varepsilon) =\sum_{k=0}^{\infty}\varepsilon^{k}  \bar{f}_i^{n,k} \quad (i=1,2,3, \ n=0,1,2,\cdots) \label{ansatz_bar_f}\\
&&V_i(\varepsilon)=  \sum_{k=0}^{\infty} \varepsilon^k V_i^k \quad (i=1,2) \label{ansatz_V}\\
&& W_i(\varepsilon)=  \sum_{k=0}^{\infty} \varepsilon^k W_i^k \quad (i=1,2) \label{ansatz_W}
\end{eqnarray} 

 
The substitution of the developments \eqref{ansatz_bar_u}-\eqref{ansatz_W}
 in \eqref{ec_uin}-\eqref{pn},  \eqref{ui0_Vi} and \eqref{sum_uik_Wi_Vi}-\eqref{sum_u3k_h}
and the identification of the terms multiplied by the same power of
$\varepsilon $, lead to a series of equations that will allow us to determine
$\bar{u}^{n,0}_i, \bar{p}^{n,-2}, \dots (i=1,2,3, \ n=0,1,2,3)$.

In this way, we identify the terms multiplied by
\begin{itemize}
    \item $\varepsilon^{-2}$ in \eqref{ec_uin} and \eqref{pn}:
    \begin{align}  &\dfrac{1}{\rho_0}  \sum_{l=1}^2\dfrac{
\partial  \bar{p}^{0,-2}}{\partial \xi_l} J_{il}^{0,0} =
 {\dfrac{ 2 \nu  }{h^2}\bar{u}_i^{2,0}} ,& (i=1,2) \label{p0-2_ui20}\\
 &\dfrac{1}{\rho_0} \sum_{l=1}^2\dfrac{
\partial  \bar{p}^{1,-2}}{\partial \xi_l} J_{il}^{0,0}  +   \bar{p}^{1,-2}    h^{-1} J_{i3}^{0,0} = \nu  
 {\dfrac{6 }{ h^2}\bar{u}_i^{3,0}}, &  (i=1,2) \label{p1-2_ui30}\\
 &\dfrac{1}{\rho_0} \sum_{l=1}^2\dfrac{
\partial  \bar{p}^{n,-2}}{\partial \xi_l} J_{il}^{0,0}  + n  \bar{p}^{n,-2}    h^{-1} J_{i3}^{0,0} = 0, & (n=2,3, \quad i=1,2) \label{pn-2_0}\\
       &\bar{p}^{n,-2}=0, & (n=1,2,3) \label{pn-2}
    \end{align}
    
Taking into account \eqref{pn-2}, we obtain from \eqref{p1-2_ui30}:
\begin{equation}
    \bar{u}_i^{3,0}=0, \quad (i=1,2) \label{ui30}
\end{equation}

    \item $\varepsilon^{-1}$ in \eqref{ec_uin} and \eqref{pn} (considering \eqref{pn-2} and \eqref{ui30}):
    \begin{align}
&{}\dfrac{1}{\rho_0}  \sum_{l=1}^2\dfrac{
\partial  \bar{p}^{0,-1}}{\partial \xi_l}  J_{il}^{0,0}=  
\dfrac{\nu}{ h} 
\bar{u}_i^{1,0} \dfrac{A^1}{A^0}+ \nu  
 {\dfrac{2}{h^2}\bar{u}_i^{2,1}}, \quad (i=1,2) \label{p0-1_ui10_ui21}\\
 &{}\dfrac{1}{ h}  \bar{u}_3^{1,0} \bar{u}_i^{1,0}  =-\dfrac{1}{\rho_0}  \sum_{l=1}^2\dfrac{
\partial  \bar{p}^{1,-1}}{\partial \xi_l}  J_{il}^{0,0}   -\dfrac{1}{\rho_0}  \sum_{l=1}^2\dfrac{
\partial  \bar{p}^{0,-2}}{\partial \xi_l}    h J_{il}^{0,1}  -\dfrac{1}{\rho_0}   \bar{p}^{1,-1}  h^{-1} J_{i3}^{0,0} \nonumber\\
&{} +  
\dfrac{2\nu}{ h} \dfrac{A^1}{A^0}
\bar{u}_i^{2,0} + \nu  
 {\dfrac{ 6}{ h^2}\bar{u}_i^{3,1}}, \quad (i=1,2) \label{ec_ns_eps-1}\\
  &{}\dfrac{1}{ h} \sum_{m=1}^2 \bar{u}_3^{m,0} (3-m) \bar{u}_i^{3-m,0}  =-\dfrac{1}{\rho_0}  \sum_{l=1}^2\dfrac{
\partial  \bar{p}^{2,-1}}{\partial \xi_l}  J_{il}^{0,0}     -\dfrac{2}{\rho_0}   \bar{p}^{2,-1}   h^{-1} J_{i3}^{0,0},\quad  (i=1,2)\\
  &{}\dfrac{1}{ h} \sum_{m=2}^3 \bar{u}_3^{m,0} (4-m) \bar{u}_i^{4-m,0}  =-\dfrac{1}{\rho_0}  \sum_{l=1}^2\dfrac{
\partial  \bar{p}^{3,-1}}{\partial \xi_l}  J_{il}^{0,0}  -\dfrac{n}{\rho_0}   \bar{p}^{3,-1}   h^{-1} J_{i3}^{0,0} ,\quad  (i=1,2, \ n=2,3)\\
 & \bar{p}^{n+1,-1}= {\dfrac{\mu}{ h} (n+2)\bar{u}_3^{n+2,0}}, \quad (n=0,1) \label{pn-1_u3n+20}\\
 & \bar{p}^{3,-1}= 0  \label{p3-1}\end{align}
 
 \item $\varepsilon^0$ in \eqref{ec_uin}-\eqref{pn},  \eqref{ui0_Vi} and \eqref{sum_uik_Wi_Vi} (keeping in mind \eqref{pn-2} and \eqref{ui30}):
 \begin{align}
&\hspace*{-0.5cm}  \dfrac{ \partial \bar{u}_i^{0,0}}{\partial t}  +
\sum_{k=1}^3 \bar{u}_k^{0,0} Q^0_{ik} -  \sum_{l=1}^2 \dfrac{
\partial \bar{u}_i^{0,0}}{\partial \xi_l}   C^0_l + \sum_{k=1}^2 \bar{u}_k^{0,0}  \sum_{l=1}^2 \left( \dfrac{
\partial \bar{u}_i^{0,0}}{\partial \xi_l} + \sum_{q=1}^3 \bar{u}_q^{0,0} H^0_{ilq} \right)   B_{lk}^{0}
\nonumber\\
&{}   =-\dfrac{1}{\rho_0}\sum_{l=1}^2\dfrac{
\partial  \bar{p}^{0,0}}{\partial \xi_l}  J_{il}^{0,0} + \nu \left\{  \sum_{m=1}^2\sum_{l=1}^2 \dfrac{
\partial^2   \bar{u}_i^{0,0}}{\partial \xi_l \partial \xi_m}
\iota_{lm}^{0}+   \sum_{k=1}^3 \sum_{l=1}^2 \dfrac{
\partial  \bar{u}_k^{0,0}}{\partial \xi_l} L_{ikl}^{0,0} +  \sum_{k=1}^3 \bar{u}_k^{0,0} S_{ik}^{0,0}
\right\}\nonumber\\
&{}+  
\dfrac{\nu}{ h} \dfrac{A^1}{A^0}
\bar{u}_i^{1,1} + \nu  
 {\dfrac{ 2}{ h^2}\bar{u}_i^{2,2}} +
 \bar{f}_i^{0,0}, \quad (i=1,2) \label{ec_ui00}\\
 &\hspace*{-0.5cm}  \dfrac{ \partial \bar{u}_i^{1,0}}{\partial t}  +
\sum_{k=1}^3 \bar{u}_k^{1,0} Q^0_{ik} -  \sum_{l=1}^2 \dfrac{
\partial \bar{u}_i^{1,0}}{\partial \xi_l}   C^0_l - \dfrac{1}{ h}  \bar{u}_i^{1,0} \left(\dfrac{\partial h}{\partial
t} + C^0_3\right) \nonumber\\
&{}+\sum_{m=0}^1 \sum_{k=1}^2 \bar{u}_k^{m,0}  \sum_{l=1}^2 \left( \dfrac{
\partial \bar{u}_i^{1-m,0}}{\partial \xi_l} + \sum_{q=1}^3 \bar{u}_q^{1-m,0} H^0_{ilq} \right)   B_{lk}^{0}
+ \sum_{k=1}^2 \bar{u}_k^{0,0}  h^{-1} B_{3k}^{0} 
 \bar{u}_i^{1,0} +\dfrac{2}{ h}  \bar{u}_3^{1,1}\bar{u}_i^{1,0}  
\nonumber\\
&{}   =-\dfrac{1}{\rho_0} \sum_{m=0}^1 \sum_{l=1}^2\dfrac{
\partial  \bar{p}^{m,m-1}}{\partial \xi_l}  h^{1-m} J_{il}^{0,1-m} 
-\dfrac{1}{\rho_0}   \bar{p}^{1,0}  h^{-1} J_{i3}^{0,0} \nonumber\\
&{} + \nu \left\{   \sum_{m=1}^2\sum_{l=1}^2 \dfrac{
\partial^2   \bar{u}_i^{1,0}}{\partial \xi_l \partial \xi_m}
\iota_{lm}^{0}+  \sum_{k=1}^3 \sum_{l=1}^2 \dfrac{
\partial  \bar{u}_k^{1,0}}{\partial \xi_l} L_{ikl}^{1,1} +  \sum_{k=1}^3 \bar{u}_k^{1,0} S_{ik}^{1,1}
\right\}\nonumber\\
&{}+  
\dfrac{2\nu}{ h} \dfrac{A^1}{A^0}
\bar{u}_i^{2,1} + \nu  
 {\dfrac{ 6 }{h^2}\bar{u}_i^{3,2}} +
 \bar{f}_i^{1,0}, \quad (i=1,2) \label{ec_ui10}\\
 &\hspace*{-0.5cm}  \dfrac{ \partial \bar{u}_i^{2,0}}{\partial t}  +
\sum_{k=1}^3 \bar{u}_k^{2,0} Q^0_{ik} -  \sum_{l=1}^2 \dfrac{
\partial \bar{u}_i^{2,0}}{\partial \xi_l}   C^0_l - \dfrac{2}{ h}  \bar{u}_i^{2,0} \left(\dfrac{\partial h}{\partial
t} + C^0_3\right) \nonumber\\
&{}+\sum_{m=0}^2 \sum_{k=1}^2 \bar{u}_k^{m,0}  \sum_{l=1}^2 \left( \dfrac{
\partial \bar{u}_i^{2-m,0}}{\partial \xi_l} + \sum_{q=1}^3 \bar{u}_q^{2-m,0} H^0_{ilq} \right)   B_{lk}^{0}
\nonumber\\
&{}+ \sum_{m=0}^{1} \sum_{k=1}^2 \bar{u}_k^{m,0}  (2-m)  h^{-1} B_{3k}^{0} 
 \bar{u}_i^{2-m,0} +\dfrac{1}{ h} \sum_{m=1}^2  (3-m)(\bar{u}_3^{m,1}\bar{u}_i^{3-m,0}  
 +\bar{u}_3^{m,0} \bar{u}_i^{3-m,1} )\nonumber\\
&{}   =-\dfrac{1}{\rho_0} \sum_{m=0}^2 \sum_{l=1}^2\dfrac{
\partial  \bar{p}^{m,m-2}}{\partial \xi_l}  h^{2-m} J_{il}^{0,2-m}  -\dfrac{1}{\rho_0} \sum_{m=1}^{2}   m  \bar{p}^{m,m-2}  h^{1-m} J_{i3}^{0,2-m} \nonumber\\
&{} + \nu \left\{ \sum_{m=1}^2\sum_{l=1}^2 \dfrac{
\partial^2   \bar{u}_i^{2,0}}{\partial \xi_l \partial \xi_m}
\iota_{lm}^{0}+  \sum_{k=1}^3 \sum_{l=1}^2 \dfrac{
\partial  \bar{u}_k^{2,0}}{\partial \xi_l} L_{ikl}^{2,2} + \sum_{k=1}^3 \bar{u}_k^{2,0} S_{ik}^{2,2}
\right\}\nonumber\\
&{}+  
\dfrac{3\nu}{ h} \dfrac{A^1}{A^0}
\bar{u}_i^{3,1} +
 \bar{f}_i^{2,0}, \quad (i=1,2)\\
 &\hspace*{-0.5cm}  \sum_{m=1}^2 \sum_{k=1}^2 \bar{u}_k^{m,0}  \sum_{l=1}^2 \left( \dfrac{
\partial \bar{u}_i^{3-m,0}}{\partial \xi_l} + \sum_{q=1}^3 \bar{u}_q^{3-m,0} H^0_{ilq} \right)   B_{lk}^{0}
\nonumber\\
&{}+ \sum_{m=1}^{2} \sum_{k=1}^2 \bar{u}_k^{m,0}  (3-m)  h^{-1} B_{3k}^{0} 
 \bar{u}_i^{3-m,0} +\dfrac{1}{ h} \sum_{m=1}^3  (4-m) (\bar{u}_3^{m,1}\bar{u}_i^{4-m,0}  
 +\bar{u}_3^{m,0} \bar{u}_i^{4-m,1} )\nonumber\\
&{}   =-\dfrac{1}{\rho_0} \sum_{m=2}^3 \sum_{l=1}^2\dfrac{
\partial  \bar{p}^{m,m-3}}{\partial \xi_l}  h^{3-m} J_{il}^{0,3-m}  
-\dfrac{1}{\rho_0} \sum_{m=2}^{3}   m  \bar{p}^{m,m-3}  h^{2-m} J_{i3}^{0,3-m}  +
 \bar{f}_i^{3,0},\nonumber\\
&{}  \quad (i=1,2)
\\
&\hspace*{-0.5cm}\bar{u}_3^{0,0}= \dfrac{\partial \vec{X}}{\partial t}
\cdot \vec{a}_3 \label{u300}\\
&\hspace*{-0.5cm}\bar{u}_3^{n+1,0}=0, \quad (n=0,1,2) \label{u3n0}\\
&\hspace*{-0.5cm}  \bar{p}^{n+1,0}= {\dfrac{\mu}{ h} (n+2)\bar{u}_3^{n+2,1}}, \quad  
 (n=0,1) \label{pn0}\\
&\hspace*{-0.5cm}  \bar{p}^{3,0}= 0 \label{p30}\\
&\hspace*{-0.5cm} \bar{u}_i^{0,0}=V_i^0 \quad (i=1,2) \label{ui00_Vi0}\\
&\hspace*{-0.5cm}  \bar{u}_i^{1,0}+\bar{u}_i^{2,0}= W_i^0 - V_i^0, \quad (i=1,2)
 \label{sum_uik0_Wi_Vi0}
\end{align}
Bearing in mind \eqref{u3n0} we yield from \eqref{pn-1_u3n+20}:
\begin{equation}
    \bar{p}^{1,-1}=\bar{p}^{2,-1}=0 \label{p1-1,p2-1}
\end{equation}
\item $\varepsilon$  in \eqref{u30m}, \eqref{u3n}, \eqref{ui0_Vi} and \eqref{sum_uik_Wi_Vi}-\eqref{sum_u3k_h} (considering \eqref{u3n0}):
\begin{align}
    &\hspace*{-0.5cm}\bar{u}_3^{0,1}= 0 \label{u301}\\
&\hspace*{-0.5cm}\bar{u}_3^{1,1}   =- h  \left[ \dfrac{1}{\sqrt{A^0}} \textrm{div}\left( \sqrt{A^0} \vec{u}^{0,0}\right) +   \bar{u}_3^{0,0}   \dfrac{A^1}{A^0}
 \right] \label{u311}\\
&\hspace*{-0.5cm}\bar{u}_3^{n+1,1}   =- \dfrac{ h}{n+1}  \left[ \dfrac{1}{\sqrt{A^0}} \textrm{div}\left( \sqrt{A^0} \vec{u}^{n,0}\right)
-  \dfrac{ n}{h}  \nabla h \cdot  \vec{u}^{n,0} 
 \right],  \quad (n=1,2) \label{u3_231}\\
&\hspace*{-0.5cm}\bar{u}_i^{0,1}=V_i^1 \quad (i=1,2) \label{ui01_Vi1}\\
&\hspace*{-0.5cm}\sum_{k=1}^3  \bar{u}_i^{k,1}= W_i^1 - V_i^1 \quad (i=1,2)
 \label{sum_uik1_Wi1_Vi1}\\
&\hspace*{-0.5cm}\sum_{k=1}^3  \bar{u}_3^{k,1}=  \dfrac{\partial
h}{\partial t}
\label{sum_u3k1_h}
\end{align} where $\vec{u}^{n,k}=(\bar{u}_1^{n,k},\bar{u}_2^{n,k})$.

\item $\varepsilon^2$ in \eqref{u30m}, \eqref{u3n} and \eqref{ui0_Vi} and \eqref{sum_uik_Wi_Vi}-\eqref{sum_u3k_h} (keeping in mind \eqref{u301}):
\begin{align}
    &\hspace*{-0.5cm}\bar{u}_3^{0,2}= 0 \label{u302}\\
&\hspace*{-0.5cm}\bar{u}_3^{1,2}   =-  \dfrac{h}{\sqrt{A^0}} \textrm{div}\left( \sqrt{A^0} \vec{u}^{0,1}\right) \label{u312} \\
&\hspace*{-0.5cm}\bar{u}_3^{2,2}   =- \dfrac{ h}{2} \sum_{m=0}^1    h^{1-m} \sum_{k=1}^2 \left[ \sum_{l=1}^2  \left(\dfrac{
\partial \bar{u}_k^{m,m}}{\partial \xi_l}   B_{lk}^{1-m}  +  \bar{u}_k^{m,m} H_{llk}^{1-m} \right)+   \bar{u}_3^{m,m}    H_{kk3}^{1-m}
 +  \dfrac{ m}{h}  \bar{u}_k^{m,m}    
B_{3k}^{1-m}  \right]\\
&\hspace*{-0.5cm}\bar{u}_3^{3,2}   =- \dfrac{ h}{3} \sum_{m=1}^2   h^{2-m} \sum_{k=1}^2 \left[ \sum_{l=1}^2  \left(\dfrac{
\partial \bar{u}_k^{m,m-1}}{\partial \xi_l}   B_{lk}^{2-m}  +  \bar{u}_k^{m,m-1} H_{llk}^{2-m} \right)+   \bar{u}_3^{m,m-1}    H_{kk3}^{2-m}\nonumber\right.\\
 &\left.{}+  \dfrac{ m}{h}  \bar{u}_k^{m,m-1}    
B_{3k}^{2-m}  \right] \label{u332}\\
&\hspace*{-0.5cm}\bar{u}_i^{0,2}=V_i^2 \quad (i=1,2) \label{ui02_Vi2}\\
&\hspace*{-0.5cm}\sum_{k=1}^3  \bar{u}_i^{k,2}= W_i^2 - V_i^2 \quad (i=1,2)
 \label{sum_uik2_Wi2_Vi2}\\
&\hspace*{-0.5cm}\sum_{k=1}^3  \bar{u}_3^{k,2}= 0
\label{sum_u3k2_h}
\end{align} 
\end{itemize}

Once the equations \eqref{p0-2_ui20}-\eqref{sum_u3k2_h} have been derived, we will proceed, in the next section, to impose the boundary conditions.

\section{Imposing boundary conditions} \label{BoundaryConditions}
\subsection{Boundary conditions leading to a lubrication problem} \label{subseccion-4-1}

Let us assume that the fluid slips at the lower
surface $(\xi _{3}=0)$, and at the upper surface $(\xi _{3}=1)$, that is, 
let us assume that the tangential velocities of the fluid at
the lower and upper surfaces are known, and, that the normal velocity of both surfaces matches the normal velocities of the fluid at the surfaces.

In this case, the terms $\bar{u}_i^0$ ($i=1,2$) are known (\eqref{ui0_Vi}),
$\bar{u}_i^{0,k}$ ($i=1,2$, $k=0,1,\dots$) are known (see \eqref{ui00_Vi0}, \eqref{ui01_Vi1}, \eqref{ui02_Vi2}), and from equations \eqref{p0-2_ui20} and \eqref{sum_uik0_Wi_Vi0},  we obtain
\begin{align}
&\bar{u}_i^{2,0} = \dfrac{h^2}{2 \mu}  \sum_{l=1}^2\dfrac{
\partial  \bar{p}^{0,-2}}{\partial \xi_l} J_{il}^{0,0},  &(i=1,2) \label{ui20}
\\
&\bar{u}_i^{1,0}=  W_i^0 - V_i^0 - \dfrac{h^2}{2 \mu}  \sum_{l=1}^2\dfrac{
\partial  \bar{p}^{0,-2}}{\partial \xi_l} J_{il}^{0,0},  &(i=1,2) \label{ui1_WVui2}
\end{align}

We can substitute  $\bar{u}^{k,1}_3$ ($k=1,2,3$) in \eqref{sum_u3k1_h}  by the expressions \eqref{u311} - \eqref{u3_231} and then, using \eqref{ui00_Vi0}, \eqref{u300}, \eqref{ui20} and \eqref{ui1_WVui2} we yield:
\begin{align}
    \dfrac{ 1}{12\mu\sqrt{A^0}} \textrm{div}\left(\sqrt{A^0} h^3\sum_{l=1}^2 \dfrac{
\partial  p^{0,-2}}{\partial \xi_l}( J^{0,0}_{1l}, J^{0,0}_{2l})\right)  - h  \left( \dfrac{\partial \vec{X}}{\partial t}
\cdot \vec{a}_3 \right)\dfrac{A^1}{A^0}\nonumber\\
{}- \dfrac{ h}{2\sqrt{A^0}} \textrm{div}(\sqrt{A^0}(\vec{W}^{0}+\vec{V}^{0})) + \dfrac{1 }{2} \nabla h \cdot (\vec{W}^{0}-\vec{V}^{0})= \dfrac{\partial
h}{\partial t} \label{ec_lub}
\end{align}

\begin{rmk}
Equation \eqref{ec_lub} is exactly the same as equation  \eqref{Reynolds_gen} (equation (99) in \cite{RodTabJMAA2021}), although the divergence term is written in a slightly different form using$$ \dfrac{M}{\sqrt{A^0}}= \sqrt{A^0}  \begin{pmatrix}
J_{11}^{0,0}& J_{12}^{0,0}\\ J_{21}^{0,0} & J_{22}^{0,0}
 \end{pmatrix} $$
As explained in \cite{RodTabJMAA2021}, \eqref{ec_lub} is a generalized version of the Reynolds equation.
\end{rmk}
 
The following result can be proved using \eqref{cambio_base_u}, \eqref{ui_pol_xi3m}, \eqref{p_pol_xi3m}, \eqref{pn-2_0}, \eqref{ui30},  \eqref{u300}, \eqref{u3n0},  \eqref{ui00_Vi0}, \eqref{ui20} and \eqref{ui1_WVui2}:
\begin{theorem} \label{Teorema-1}
If we assume that there exist asymptotic expansions  \eqref{ansatz_bar_u}-\eqref{ansatz_W}, that \eqref{f_pol_xi3} holds
 and that the tangential and normal velocities are known on the 
bound surfaces, then the solution of model \eqref{ui_pol_xi3m}-\eqref{pn}  verifies 
\begin{align}
u_k(\varepsilon) &= V_k^0 + \xi_3 (W_k^0 - V_k^0)  +\dfrac{h^2}{2 \mu}  \sum_{l=1}^2\dfrac{
\partial  \bar{p}^{0,-2}}{\partial \xi_l} J_{kl}^{0,0}(\xi_3^2-\xi_3) + O(\varepsilon), 
\quad (k=1,2) \label{eq-1-T-1} \\
u_3(\varepsilon) &= \dfrac{\partial \vec{X}}{\partial t}
\cdot \vec{a}_3 + O(\varepsilon), \\
p(\varepsilon) &= \varepsilon^{-2} \bar{p}^{0,-2} + O(\varepsilon^{-1}), \label{eq-3-T-1}
\end{align}
where $\bar{p}^{0,-2}$is the solution of the equation \eqref{ec_lub}. Thus, 
the velocity, $\vec{u}^{\varepsilon}$, and the pressure, $p^{\varepsilon}$, defined in the original domain, satisfy 
\begin{align}
u_i^{\varepsilon} &= \sum_{k=1}^2 u_k(\varepsilon) \dfrac{\partial x_i}{\partial \xi_k} + u_3(\varepsilon) N_i, \quad (i=1,2,3)    \\
&=\sum_{k=1}^2 \left[V_k^0 + \xi_3 (W_k^0 - V_k^0)  +\dfrac{h^2}{2 \mu}  \sum_{l=1}^2\dfrac{
\partial  \bar{p}^{0,-2}}{\partial \xi_l} J_{il}^{0,0}(\xi_3^2-\xi_3) \right] \dfrac{\partial x_i}{\partial \xi_k} + \left( \dfrac{\partial \vec{X}}{\partial t}
\cdot \vec{a}_3\right) N_i + O(\varepsilon),\\
p^{\varepsilon}&=\varepsilon^{-2} \bar{p}^{0,-2} + O(\varepsilon^{-1}).
\end{align}
\end{theorem}

\begin{rmk}\label{rmkth1}
Following the steps of \cite{CR2}, it would be possible to show that the norm of the difference between the exact solution of \eqref{ui_pol_xi3m}-\eqref{pn} and its asymptotic approximation is small. In this way, error estimates could be obtained not just formally and,  Theorem~\ref{Teorema-1} could be mathematically proved. However, the adaptation of the method used in \cite{CR2} to prove Theorem~\ref{Teorema-1} is not trivial and it would be very laborious.
\end{rmk}

\subsection{Boundary conditions leading to a thin fluid layer problem} \label{subseccion_sw}
 
Now, instead
of considering that the tangential and normal velocities are known on the
upper and lower surfaces,
we assume that the normal component of the traction on $\xi _{3}=0$ and
on $\xi _{3}=1$ are known pressures, and that the tangential component
of the traction on these surfaces are friction forces depending on the
value of the velocities on $\partial D$. Therefore, we assume that
\begin{align}
&\vec{T}^{\varepsilon}\cdot
\vec{n}^{\varepsilon}_0 = (\sigma^{\varepsilon} \vec{n}^{\varepsilon}_0)\cdot
\vec{n}^{\varepsilon}_0=-\pi^{\varepsilon}_0 &\textrm{ on } \xi_3=0& \label{Tn0_xi3_0}\\
&\vec{T}^{\varepsilon}\cdot
\vec{n}^{\varepsilon}_1 = (\sigma^{\varepsilon} \vec{n}^{\varepsilon}_1)\cdot
\vec{n}^{\varepsilon}_1=-\pi^{\varepsilon}_1 &\textrm{ on } \xi_3=1&
\label{Tn1_xi3_1}
\\
&\vec{T}^{\varepsilon}\cdot
\vec{a}_i= (\sigma^{\varepsilon}  \vec{n}^{\varepsilon}_0)\cdot
\vec{a}_i=-\vec{f}^{\varepsilon}_{R_0}\cdot
\vec{a}_i &\textrm{ on } \xi_3=0, &\quad (i=1,2) \label{Tai_xi3_0}\\
&\vec{T}^{\varepsilon}\cdot
\vec{v}_i^{\varepsilon}= (\sigma^{\varepsilon}   \vec{n}^{\varepsilon}_1)\cdot
\vec{v}_i^{\varepsilon}=-\vec{f}^{\varepsilon}_{R_1}\cdot
\vec{v}_i^{\varepsilon} &\textrm{ on } \xi_3=1, &\quad (i=1,2)
\label{Tvi_xi3_1}
\end{align} where $\vec{T}^{\varepsilon}$ is the traction vector and $\sigma^{\varepsilon}$ is the stress tensor given by
\begin{align}
    \sigma^{\varepsilon}_{ij}&= -p^{\varepsilon}\delta_{ij}+\mu \left(
\dfrac{\partial u_i^{\varepsilon}}{\partial x^{\varepsilon}_j} +
\dfrac{\partial u_j^{\varepsilon}}{\partial
x^{\varepsilon}_i}\right)\nonumber\\
&= \sum_{n=0}^3 \xi_3^n \left[ -\sum_{k=-2}^{\infty} \varepsilon^k \bar{p}^{n,k} \delta_{ij}+\mu \sum_{k=0}^{\infty} \varepsilon^k\sum_{m=1}^3 \sum_{l=1}^2 \left(
\dfrac{\partial( \bar{u}_m^{n,k}  a_{mi})}{\partial \xi_l} \left(\sum_{r=0}^{\infty}(\varepsilon \xi_3 h)^r(\alpha_l^r a_{1j}+\beta_l^r a_{2j})\right) \right. \right.\nonumber\\
&+\left. \left.
\dfrac{\partial (\bar{u}_m^{n,k} a_{mj})}{\partial
\xi_l}  \left(\sum_{r=0}^{\infty}(\varepsilon \xi_3 h)^r(\alpha_l^r a_{1i}+\beta_l^r a_{2i})\right) \right)\right]\nonumber\\
&+ \mu \sum_{k=0}^{\infty} \varepsilon^k \sum_{n=1}^3 \sum_{m=1}^3 n \xi_3^{n-1}  \bar{u}_m^{n,k} \left[\sum_{r=0}^{\infty}\varepsilon^r \xi_3^{r+1} h^{r-1}  
 \left(\alpha_3^r(  a_{mi} a_{1j}+ a_{mj} a_{1i}) \right.\right.\nonumber\\
&+\left.\left.\beta_3^nr(a_{mi} a_{2j}+ a_{mj} a_{2i}) \right) + \dfrac{1}{\varepsilon h} ( a_{mi}  a_{3j} +  a_{mj}a_{3i})  \right] ,  \quad (i,j=1,2,3) \label{sigmaij_des}
\end{align}
and vectors $\vec{n}^{\varepsilon}_0$, $\vec{n}^{\varepsilon}_1$ are, respectively, the outward unit normal vectors to the lower and the upper surfaces, that is 
\begin{align}
&\vec{n}^{\varepsilon}_0=s_0 \vec{a}_3 \label{n0}\\
&\vec{n}^{\varepsilon}_1=-s_0\dfrac{\vec{v}^{\varepsilon}_3}{\|\vec{v}^{\varepsilon}_3\|}
\label{n1} 
\end{align}
where
\begin{equation}
    s_0=-1 \textrm{ or } s_0=1 
\end{equation}
is fixed ($\vec{n}^{\varepsilon}_0 = \vec{a}_3$ or $\vec{n}^{\varepsilon}_0 = - \vec{a}_3$, depending on the orientation of the parametrization $\vec{X}$), and 
\begin{align}
\vec{v}^{\varepsilon}_1&=\vec{a}_1+\varepsilon
\left(\dfrac{\partial h}{\partial \xi_1} \vec{a}_3 + h
\dfrac{\partial \vec{a}_3}{\partial \xi_1}\right) \label{v1}\\
\vec{v}^{\varepsilon}_2&=\vec{a}_2+\varepsilon
\left(\dfrac{\partial h}{\partial \xi_2} \vec{a}_3 + h
\dfrac{\partial \vec{a}_3}{\partial \xi_2}\right) \label{v2}\\
\vec{v}^{\varepsilon}_3&=\vec{v}^{\varepsilon}_1 \times \vec{v}^{\varepsilon}_2 =\vec{a}_1 \times \vec{a}_2 +\varepsilon
\left[\dfrac{\partial h}{\partial \xi_2} (\vec{a}_1 \times
\vec{a}_3) + h \left(\vec{a}_1 \times \dfrac{\partial
\vec{a}_3}{\partial \xi_2} + \dfrac{\partial
\vec{a}_3}{\partial \xi_1} \times
\vec{a}_2 \right) + \dfrac{\partial h}{\partial
\xi_1} (\vec{a}_3\times \vec{a}_2) \right]\nonumber\\
&+\varepsilon^2 \left[ \left(\dfrac{\partial h}{\partial \xi_1}
\vec{a}_3 + h \dfrac{\partial \vec{a}_3}{\partial
\xi_1}\right)\times\left(\dfrac{\partial h}{\partial \xi_2}
\vec{a}_3 + h \dfrac{\partial \vec{a}_3}{\partial \xi_2}\right)
\right]\label{v3}\\
 \|\vec{v}^{\varepsilon}_3\|&=\|\vec{a}_1 \times
\vec{a}_2\| + \varepsilon h  \left[ \vec{a}_3 \cdot \left(\vec{a}_1
\times \dfrac{\partial \vec{a}_3}{\partial \xi_2}\right) +\vec{a}_3
\cdot \left(\dfrac{\partial \vec{a}_3}{\partial \xi_1} \times
\vec{a}_2\right)\right] + O(\varepsilon^2)\label{mod_v3}
\end{align}

Typically, the friction force is of the form
\begin{equation}
    \vec{f}^{\varepsilon}_{R\alpha} = \rho_0 C_R^\varepsilon \|
\vec{{u}}^\varepsilon \| \vec{{u}}^\varepsilon \textrm{ on } \xi_3=\alpha, \quad (\alpha=0,1)\label{fR}
\end{equation}
where $C_R^\varepsilon$ is a small constant. Let us assume that it is of order $\varepsilon$ (see \cite{RTV1} or \cite{Weiyan}), that is, 
\begin{equation}
    C_R^{\varepsilon} =\varepsilon C^1_R\label{CR}
\end{equation}

If we assume that the pressures and the friction forces on the upper and lower surfaces admit a development in powers of $\varepsilon$ too:
\begin{align}
\pi_i(\varepsilon)&=\sum_{r=0}^{\infty} \varepsilon^r \pi_i^r, \quad (i=0,1) \label{ansatz_pi}\\
 \vec{f}(\varepsilon)_{R_\alpha}&=\sum_{k=1}^{\infty} \varepsilon^k  \vec{f}^k_{R_\alpha}, \quad (\alpha=0,1) \label{ansatz_fR}
\end{align}
condition \eqref{Tn0_xi3_0} can now be written (using
\eqref{sigmaij_des}, \eqref{n0}) as:
\begin{equation}
(\sigma_{ij}a_{3j})a_{3i} = -\bar{p}^0 + \dfrac{2\mu }{\varepsilon h}  \bar{u}_3^1=-\pi_0 \label{Tn0_xi3_0s}
\end{equation}

Next, we identify  the terms multiplied by the same power of $\varepsilon$, we obtain (for $\varepsilon^{-2}$, $\varepsilon^{-1}$ and $\varepsilon^0$) taking into account  \eqref{u300}, \eqref{u3n0}  and \eqref{u311}:
\begin{align}
&\bar{p}^{0,-2}=\bar{p}^{0,-1}=0 \label{p0-2-1_sw}\\
&\bar{p}^{0,0}=\pi_0^0-\dfrac{2\mu}{\sqrt{A^0}}\textrm{div}(\sqrt{A^0}\vec{V}^{0})- 2\mu \left( \dfrac{\partial \vec{X}}{\partial t} \cdot \vec{a}_3\right) \dfrac{A^1}{A^0}
 \label{p00_sw}
\end{align}
 
From \eqref{p0-2_ui20}, \eqref{p0-2-1_sw}, \eqref{u3_231}, \eqref{pn0}, \eqref{ec_ns_eps-1}, \eqref{sum_uik0_Wi_Vi0}, \eqref{p0-1_ui10_ui21} and bearing in mind \eqref{u3n0}, \eqref{p1-1,p2-1}, \eqref{p0-2-1_sw}  we get:
\begin{align}
\bar{u}_i^{2,0}&=0, \quad (i=1,2) \label{ui20_sw}
\\
\bar{u}_3^{3,1}&=0 \label{u331_sw}\\
\bar{p}^{2,0}&=0 \label{p20_sw}\\
\bar{u}_i^{3,1}&= 0  \quad (i=1,2) \label{ui31_sw}\\
\bar{u}_i^{1,0}&=W_i^0-V_i^0, \quad (i=1,2) \label{ui10_sw}\\
\bar{u}_i^{2,1}&=-
 \dfrac{h}{ 2} \dfrac{A^1}{A^0}
(W_i^0-V_i^0). \quad (i=1,2)   
 \label{ui21_ui10}
\end{align}

Boundary condition \eqref{Tn1_xi3_1} can be written (using
 \eqref{n1}) as follows:
\begin{eqnarray}
&&\left(\sigma_{ij}^{\varepsilon}
{v}^{\varepsilon}_{3j}\right)\cdot
{v}^{\varepsilon}_{3i}=-\pi^{\varepsilon}_1\|\vec{v}^{\varepsilon}_3\|^2
\textrm{ on } \xi_3=1 \label{Tn1_xi3_1s}
\end{eqnarray}
and, using
\eqref{sigmaij_des},   \eqref{v3}, \eqref{mod_v3} and \eqref{ansatz_pi} to substitute $\sigma_{ij}$, vector $\vec{v}_3$, its module and $\pi^{\varepsilon}_1$ into the above condition, we identify the terms multiplied by $\varepsilon^0$ (the terms multiplied by $\varepsilon^{-2}$ and $\varepsilon^{-1}$ trivially vanish):
\begin{align}
&
- \sum_{n=0}^3 \bar{p}^{n,0} \|\vec{a}_1 \times \vec{a}_2\|^2 + \dfrac{ 2\mu}{ h}  \sum_{n=1}^3  n  \bar{u}_3^{n,1} \|\vec{a}_1 \times \vec{a}_2\|^2\nonumber\\
&{}+  \dfrac{2 \mu}{ h}  \sum_{n=1}^3 \sum_{m=1}^2 n  \bar{u}_m^{n,0}  \|\vec{a}_1 \times \vec{a}_2\| \vec{a}_{m} \cdot
\left[\dfrac{\partial h}{\partial \xi_2} (\vec{a}_1 \times
\vec{a}_3) + h \left(\vec{a}_1 \times \dfrac{\partial
\vec{a}_3}{\partial \xi_2} + \dfrac{\partial
\vec{a}_3}{\partial \xi_1} \times
\vec{a}_2\right)+ \dfrac{\partial h}{\partial
\xi_1} (\vec{a}_3\times \vec{a}_2)\right]\nonumber\\
&=- \pi_1^0\|\vec{a}_1 \times\vec{a}_2\|^2 \label{Tni1_xi3_1_eps0}
\end{align}

Taking into account \eqref{p20_sw}, \eqref{p30}, \eqref{ui30}, \eqref{ui20_sw}, \eqref{u331_sw}, \eqref{pn0}, \eqref{u3_231}, \eqref{p00_sw}, \eqref{ui10_sw}  and, that
\begin{align}
&(\vec{a}_1 \times \vec{a}_3)\cdot\vec{a}_{2}  = (\vec{a}_3\times
\vec{a}_2) \cdot \vec{a}_1 =- \| \vec{a}_{1}\times \vec{a}_2\| \label{a1xa3_a2}\\
& \left(\vec{a}_1 \times \dfrac{\partial \vec{a}_3}{\partial
\xi_2}\right)\cdot\vec{a}_{2} = \left(\dfrac{\partial \vec{a}_3}{\partial \xi_1} \times
\vec{a}_2\right)\cdot \vec{a}_1 = 0 \label{a12_x_da3dxi21_a21}
\end{align}
equation \eqref{Tni1_xi3_1_eps0} can be written as follows:
\begin{align}
&{}\pi_0^0
 + \dfrac{ \mu}{ h}   \left(
   \dfrac{ h}{\sqrt{A^0}} \textrm{div}(\sqrt{A^0}(\vec{W}^{0}-\vec{V}^0)) + \nabla h \cdot (\vec{W}^{0}-\vec{V}^0)\right)
   = \pi_1^0\label{Tni1_xi3_1_eps0s}
\end{align}

Boundary conditions \eqref{Tai_xi3_0} on $\xi_3=0$ have been re-written  using
\eqref{sigmaij_des}, \eqref{n0} and \eqref{ansatz_fR}, then we identify the terms multiplied by each power of $\varepsilon$:
\begin{itemize}
\item $\varepsilon^{-1}$: 
\begin{equation}
    \bar{u}_i^{1,0}=0, \quad (i=1,2) \label{ui10_sw0} 
\end{equation}
and from \eqref{ui10_sw}:
\begin{equation}
W_i^0=V_i^0, \quad (i=1,2) \label{W_V_sw}
\end{equation}
\item $\varepsilon^0$ (keeping in mind \eqref{u300}, \eqref{ui00_Vi0}):
\begin{equation}
\dfrac{\partial }{\partial
\xi_i} \left(  \dfrac{\partial \vec{X}}{\partial t}
\cdot \vec{a}_3\right)+ \sum_{m=1}^2 V_m^{0} \left(
\dfrac{\partial  \vec{a}_{m}}{\partial
\xi_i} \cdot \vec{a}_{3}\right) + \dfrac{1}{ h}\sum_{m=1}^2   \bar{u}_m^{1,1}   \vec{a}_{m}\cdot \vec{a}_{i} =0, \quad (i=1,2) \label{cc_xi3_0_fR0_eps0}
\end{equation}
\item $\varepsilon$ (considering \eqref{u301}, \eqref{ui01_Vi1}):
\begin{equation}
s_0  \mu \left[  
\sum_{m=1}^2  V_m^{1} \left(
\dfrac{\partial  \vec{a}_{m}}{\partial
\xi_i} \cdot \vec{a}_{3}\right)+  \dfrac{1}{ h} \sum_{m=1}^2   \bar{u}_m^{1,2}   \vec{a}_{m}\cdot \vec{a}_{i}  \right]=-  \vec{f}^1_{R_0}\cdot \vec{a}_{i}. \quad (i=1,2) \label{cc_xi3_0_fR0_epsk}
\end{equation}
\end{itemize}

If we sum equation \eqref{cc_xi3_0_fR0_eps0} ($i=1$)  multiplied   by $\alpha^0_1$ and equation \eqref{cc_xi3_0_fR0_eps0} ($i=2$) multiplied by $\beta^0_1$  (and we repeat the process analogously multiplying by $\alpha^0_2$ and $\beta^0_2$) we obtain:
\begin{equation}
\bar{u}_i^{1,1} =-h\sum_{l=1}^2\left( J^{0,0}_{il}
\dfrac{\partial }{\partial
\xi_l}\left( \dfrac{\partial \vec{X}}{\partial t}
\cdot \vec{a}_3 \right) +  V_l^{0}D^0_{il} \right), \quad (i=1,2) 
 \label{ui11_sw}
\end{equation} 
Coefficients $D^0_{il}$ are defined in \eqref{D}.

We do the same operations from \eqref{cc_xi3_0_fR0_epsk} to get:
\begin{equation}
    \bar{u}_i^{1,2}   =- h \sum_{m=1}^2 \left( \dfrac{s_0}{\mu}J^{0,0}_{im} (\vec{f}^1_{R_0}\cdot \vec{a}_{m}) +   V_m^{1} D^0_{im}\right), \quad (i=1,2)\label{ui1k_sw}
\end{equation}

From \eqref{u3_231}, \eqref{pn0}, \eqref{ui21_ui10}, \eqref{u332}, \eqref{ec_ui10}, taking into account  \eqref{ui10_sw0}, \eqref{W_V_sw}, \eqref{p0-2-1_sw}, \eqref{u3n0},  we have:
\begin{align}
\bar{u}_3^{2,1}&=0\label{u321_sw}\\
 \bar{p}^{1,0}&= 0 \label{p10_sw}\\
  \bar{u}_i^{2,1}&=0, \quad (i=1,2) \label{ui21_sw}\\
    \bar{u}_3^{3,2}&=0 \label{u332_sw}
\\
\bar{u}_i^{3,2}& =- \dfrac{  h^2 }{6\nu}\bar{f}_i^{1,0}, \quad (i=1,2) \label{ui32_fi1}
\end{align}

Next, we yield from \eqref{sum_u3k1_h}, \eqref{u311}, \eqref{u300}, \eqref{ui00_Vi0}, \eqref{u321_sw} and \eqref{u331_sw}:
\begin{align}
 \dfrac{\partial
h}{\partial t} &=\bar{u}_3^{1,1}\label{u311_sw}\\
&=- { h \left( \dfrac{\partial \vec{X}}{\partial t}
\cdot \vec{a}_3 \right)\dfrac{A^1}{A^0}} - \dfrac{ h}{\sqrt{A^0}} \textrm{div}(\sqrt{A^0}\vec{V}^0)   \label{ec_h_V_sw}
\end{align}
and, now,  \eqref{p00_sw} reads (using \eqref{ec_h_V_sw})
\begin{equation}
\bar{p}^{0,0}=\pi_0^0 + \dfrac{2\mu }{ h}  \dfrac{\partial h}{\partial t} \label{p00_sw_h}
\end{equation}

Boundary conditions \eqref{Tvi_xi3_1} can be rewritten taking into account \eqref{sigmaij_des}, \eqref{n1}, \eqref{v1}-\eqref{mod_v3} and \eqref{ansatz_fR}. Then we identify the terms multiplied by the each power of $\varepsilon$:
\begin{itemize}
    \item $\varepsilon^{0}$ (considering \eqref{ui30}, \eqref{ui20_sw}, \eqref{ui10_sw0}, \eqref{u3n0},  \eqref{ui31_sw}, \eqref{ui21_sw}) we re-obtain  \eqref{cc_xi3_0_fR0_eps0}
\item $\varepsilon$ (using \eqref{a1xa3_a2}-\eqref{a12_x_da3dxi21_a21}, bearing in mind \eqref{u301}, \eqref{u331_sw}, \eqref{u321_sw}, \eqref{ui31_sw},  \eqref{u3n0}, \eqref{ui30}, \eqref{ui20_sw}, \eqref{ui10_sw0}, \eqref{ui21_sw} and dividing by $\|\vec{a}_1 \times \vec{a}_2\|$):
\begin{align}
&\hspace*{-0.5cm} \mu\left\{\dfrac{\partial \bar{u}_3^{1,1}}{\partial
\xi_i}  +   
\sum_{m=1}^2 \bar{u}_m^{0,1} \left(
\dfrac{\partial  \vec{a}_{m}}{\partial
\xi_i}  \cdot \vec{a}_{3}\right) 
   + \dfrac{1}{ h}\sum_{n=1}^3 \sum_{m=1}^2 n  \bar{u}_m^{n,2}   (\vec{a}_{m}\cdot\vec{a}_i)   + \dfrac{1}{ h}\dfrac{\partial h}{\partial \xi_i}  \bar{u}_3^{1,1}   \right. \nonumber\\
&\left.{}-  \sum_{l=1}^2  \left(\alpha_l^0 \dfrac{\partial h}{\partial\xi_1}+\beta_l^0 \dfrac{\partial h}{\partial\xi_2}\right) \left[\sum_{m=1}^2 
\dfrac{\partial \bar{u}_m^{0,0} }{\partial \xi_l}\left( \vec{a}_m \cdot \vec{a}_{i}  \right)+ \sum_{m=1}^3 \bar{u}_m^{0,0} \left(\dfrac{\partial \vec{a}_{m}}{\partial \xi_l} \cdot  \vec{a}_i\right)  \right]\right.\nonumber\\
&{}
-\sum_{m=1}^2
\dfrac{\partial \bar{u}_m^{0,0}}{\partial
\xi_i} \dfrac{\partial h}{\partial \xi_m} +\dfrac{ h}{\sqrt{A^0}} \dfrac{\partial \bar{u}_3^{0,0}}{\partial
\xi_i} I + \dfrac{1}{\sqrt{A^0}} \sum_{m=1}^3  \bar{u}_m^{0,0} 
\left( \dfrac{\partial  \vec{a}_{m}}{\partial
\xi_i} \cdot \vec{\eta}(h)\right)\nonumber\\
& \left.{}+   \dfrac{1}{\sqrt{A^0}} \sum_{m=1}^2  \bar{u}_m^{1,1}( \vec{a}_{m}\cdot  \vec{a}_{s})I\right\}=s_0 \left( \vec{f}^1_{R_1}\cdot
\vec{a}_i\right) \quad (i=1,2) \label{Tvi_xi3_1_epsb}
\end{align} where coefficients $I$ and $\vec{\eta}(h)$ are defined in \eqref{I} and \eqref{eta} respectively.
\end{itemize}

Now, we replace in \eqref{Tvi_xi3_1_epsb} the terms $\bar{u}_3^{1,1}$, $\bar{u}_3^{0,0}$, $\bar{u}_i^{0,0}$, $\bar{u}_i^{0,1}$, $\bar{u}_i^{1,1}$ and $\bar{u}_i^{1,2}$ ($i=1,2$) by the expressions obtained in \eqref{u311_sw}, \eqref{u300}, \eqref{ui00_Vi0}, \eqref{ui01_Vi1}, \eqref{ui11_sw} and \eqref{ui1k_sw}  respectively, and taking into account \eqref{ui32_fi1}, we can write:
\begin{align}
&\hspace*{-0.5cm} \dfrac{\partial^2 h}{\partial t\partial
\xi_i}  + \dfrac{1}{ h}\dfrac{\partial h}{\partial \xi_i}  \dfrac{\partial h}{\partial t}
+ \dfrac{2}{ h} \sum_{m=1}^2   \bar{u}_m^{2,2}   (\vec{a}_{m}\cdot\vec{a}_i)   \nonumber\\
&\left.{}-  \sum_{l=1}^2  \left(\alpha_l^0 \dfrac{\partial h}{\partial\xi_1}+\beta_l^0 \dfrac{\partial h}{\partial\xi_2}\right) \left[\sum_{m=1}^2 
\dfrac{\partial V_m^{0} }{\partial \xi_l} \left( \vec{a}_m \cdot \vec{a}_{i}  \right)+ \sum_{m=1}^3 V_m^{0}\left(\dfrac{\partial \vec{a}_{m}}{\partial \xi_l} \cdot  \vec{a}_i\right)  \right]\right.\nonumber\\
&{} -\sum_{m=1}^2
\dfrac{\partial V_m^{0}}{\partial
\xi_i} \dfrac{\partial h}{\partial \xi_m} - \dfrac{1}{\sqrt{A^0}} \sum_{m=1}^2 V_m^{0}\left[
\left( \dfrac{\partial  \vec{a}_{m}}{\partial
\xi_i} \cdot \vec{\eta}(h)\right) -  h I \left( \vec{a}_{3}\cdot  \dfrac{\partial \vec{a}_m}{\partial\xi_i} \right)\right]\nonumber\\
& {}
+  \left( \dfrac{\partial \vec{X}}{\partial t}
\cdot \vec{a}_3 \right)\left[ H^0_{3i3} + \dfrac{1}{\sqrt{A^0}}  
\left( \dfrac{\partial  \vec{a}_{3}}{\partial
\xi_i} \cdot \vec{\eta}(h)\right)\right] \nonumber\\
&{}= \dfrac{s_0}{\mu} \left( \vec{f}^1_{R_1}+ \vec{f}^1_{R_0} \right)\cdot
\vec{a}_i  +\dfrac{h}{2\nu}  \sum_{m=1}^2   \bar{f}_m^{1,0}  (\vec{a}_{m}\cdot\vec{a}_i)     \quad (i=1,2) \label{Tvi_xi3_1_epse}
\end{align}

Next, we multiply \eqref{Tvi_xi3_1_epse}  by $J^{0,0}_{ji}$ for $j=1,2$ and we sum in $i=1,2$. In this way we are able to infer the terms 
$\bar{u}_i^{2,2}$
 \begin{align}
&\hspace*{-0.5cm} \dfrac{2}{ h}   \bar{u}_i^{2,2}  
=-  \sum_{l=1}^2 
\dfrac{\partial V_i^{0} }{\partial \xi_l} J^{0,0}_{3l}
- \sum_{m=1}^2 V_m^{0}\left(\sum_{l=1}^2 J^{0,0}_{3l} H^0_{ilm} + \dfrac{ h I}{\sqrt{A^0}} 
   D^0_{im}\right)   \nonumber\\
&{}-\sum_{j=1}^2 J^{0,0}_{ij}\left(\dfrac{\partial^2 h}{\partial t\partial
\xi_j}  + \dfrac{1}{ h}\dfrac{\partial h}{\partial \xi_j}  \dfrac{\partial h}{\partial t} -\sum_{m=1}^2
\dfrac{\partial V_m^{0}}{\partial
\xi_j} \dfrac{\partial h}{\partial \xi_m} -  \dfrac{1}{\sqrt{A^0}} \sum_{m=1}^2 V_m^{0} \left( \dfrac{\partial  \vec{a}_{m}}{\partial
\xi_j} \cdot \vec{\eta}(h)\right)\right)   \nonumber\\
& {}
-  \left( \dfrac{\partial \vec{X}}{\partial t}
\cdot \vec{a}_3 \right)\sum_{j=1}^2 J^{0,0}_{ij}\left[ H^0_{3j3} + \dfrac{1}{\sqrt{A^0}}  
\left( \dfrac{\partial  \vec{a}_{3}}{\partial
\xi_j} \cdot \vec{\eta}(h)\right)\right]  \nonumber\\
&{}+\dfrac{s_0}{\mu} \sum_{j=1}^2 J^{0,0}_{ij}\left( \vec{f}^1_{R_1}+ \vec{f}^1_{R_0}\right) \cdot
\vec{a}_j  +\dfrac{  h}{2 \nu}   \bar{f}_i^{1,0}, \quad (i=1,2)    \label{ui22_sw}
\end{align}

We can rewrite equations \eqref{ec_ui00} considering \eqref{ui00_Vi0}, \eqref{u300}, \eqref{p00_sw_h}, \eqref{ui11_sw} and \eqref{ui22_sw}  to substitute $\bar{u}_i^{0,0}$ ($i=1,2,3$), $\bar{p}^{0,0}$, $\bar{u}_i^{1,1}$ ($i=1,2$) and $\bar{u}_i^{2,2}$ ($i=1,2$) respectively by the expressions obtained previously:
\begin{eqnarray}
&&\hspace*{-0.5cm} \dfrac{ \partial V_i^{0}}{\partial t}+\sum_{l=1}^2  \dfrac{ \partial V_i^{0}}{\partial \xi_l} (V_l^{0} -C^0_l)+ \sum_{k=1}^2
 V_k^{0}  \left(  R^0_{ik} + \sum_{l=1}^2 V_l^{0} H^0_{ilk} \right)
 \nonumber\\
&&{}
 =-\dfrac{1}{\rho_0} \sum_{l=1}^2  J^{0,0}_{il}\dfrac{
\partial \pi_0^0 }{\partial \xi_l}  + \nu \left[ \sum_{m=1}^2 \sum_{l=1}^2 \dfrac{
\partial^2 V_i^{0}}{\partial \xi_l \partial \xi_m} J^{0,0}_{lm}+ \sum_{k=1}^2 \sum_{l=1}^2 \dfrac{
\partial V_k^{0}}{ \partial \xi_l} \bar{L}_{ikl}^{0,0}  +  \sum_{k=1}^2  V^{0}_k \bar{S}^{0,0}_{ik}+\kappa^0_i\right] \nonumber\\
&&{} + \bar{F}^0_i - 
 \left( \dfrac{\partial \vec{X}}{\partial t}
\cdot \vec{a}_3 \right)  Q^0_{i3}  \quad (i=1,2) \label{ec_ui00_sw}\end{eqnarray}
where coefficients $\bar{L}^{0,0}_{ikl}$, $R^0_{ik}$, $\bar{S}^{0,0}_{ik}$, $\kappa^0_i$ and $\bar{F}^0_i$ are defined in \eqref{Lbar}, \eqref{R},  \eqref{Sbar}, \eqref{kappa} and \eqref{Fbar2} respectively.

\begin{rmk}
Equation \eqref{ec_ui00_sw} is exactly  the same as equation \eqref{ec_Vi0-v2} (see \eqref{F}-\eqref{Fbar3}).
\end{rmk}

Taking into account \eqref{cambio_base_u}, \eqref{ui_pol_xi3}, \eqref{p_pol_xi3}, \eqref{pn-2_0}, \eqref{ui30},  \eqref{u300}, \eqref{u3n0},  \eqref{ui00_Vi0}, \eqref{p0-2-1_sw}, \eqref{ui20_sw}, \eqref{ui10_sw0}, \eqref{p00_sw_h}, \eqref{p10_sw}, \eqref{p20_sw} and \eqref{p30}, we can prove the following result: 
\begin{theorem} \label{Teorema-2}
If we assume that there exist asymptotic expansions \eqref{ansatz_bar_u}-\eqref{ansatz_W} and \eqref{ansatz_pi}-\eqref{ansatz_fR}, and that the hypothesis \eqref{f_pol_xi3} and \eqref{CR}, and the boundary conditions \eqref{Tn0_xi3_0}-\eqref{Tvi_xi3_1} hold, then the solution of model \eqref{ui_pol_xi3m}-\eqref{pn} verifies
\begin{align}
u_k(\varepsilon) &= V_k^0 + O(\varepsilon), \quad (k=1,2) \label{eq-1-T-2} \\
u_3(\varepsilon) &= \dfrac{\partial \vec{X}}{\partial t}
\cdot \vec{a}_3 + O(\varepsilon), \\
p(\varepsilon) &= \pi_0^0 + \dfrac{2\mu}{h}\dfrac{\partial h}{\partial t} + O(\varepsilon), \label{eq-3-T-2}
\end{align}
where $V_k^0$ ($k=1,2$) are the solutions of the equations \eqref{ec_ui00_sw}. Thus, 
the velocity, $\vec{u}^{\varepsilon}$, and the pressure, $p^{\varepsilon}$, defined in the original domain, satisfy
\begin{align}
u_i^{\varepsilon} &=\sum_{k=1}^2 V_k^0 \dfrac{\partial x_i}{\partial \xi_k} + \left( \dfrac{\partial \vec{X}}{\partial t}
\cdot \vec{a}_3\right) N_i + O(\varepsilon), \quad (i=1,2,3)    \\
p^{\varepsilon}&= \pi_0^0 + \dfrac{2\mu }{ h^{\varepsilon}}  \dfrac{\partial h^{\varepsilon}}{\partial t} + O(\varepsilon).
\end{align}
\end{theorem}
\begin{rmk} \label{rmkth2}
As stated in Remark \ref{rmkth1}, the steps of \cite{CR2} could be followed to prove \eqref{eq-1-T-2}-\eqref{eq-3-T-2} obtaining error estimates.
\end{rmk}

\section{Conclusions}\label{sec-conclusions}

In this article, we propose a two-dimensional flow model of a viscous fluid between two very close moving surfaces and we show (using a formal asymptotic expansion of the solution) that its asymptotic behavior, when the distance between the two surfaces tends to zero, is the same as that previously obtained in \cite{RodTabJMAA2021} for the Navier-Stokes equations. In fact, we have justified that, under the assumptions about the boundary conditions made in subsection \ref{subseccion-4-1}, the solution of the new model approaches the solution of model \eqref{ec_lub} as $\varepsilon$ tends to zero, just as in the previous work \cite{RodTabJMAA2021}, where we showed that the solution of the Navier-Stokes equations
approaches the solution of \eqref{Reynolds_gen}.  And, we have also seen that, under the assumptions about the boundary conditions shown in subsection \ref{subseccion_sw}, the solution of the the new model tends to the solution of \eqref{ec_ui00_sw}, as it happened in our prior article \cite{RodTabJMAA2021} with the solution of the Navier-Stokes equations (see \eqref{ec_Vi0-v2}).

As we have already pointed out in Remarks \ref{rmkth1} and \ref{rmkth2}, the justification of Theorems \ref{Teorema-1} and \ref{Teorema-2} is based on the formal asymptotic expansion of the solution of   model \eqref{ui_pol_xi3m}-\eqref{pn}. This could be seen as an engineering approach to the justification of the newly derived model. We have also commented, in Remarks \ref{rmkth1} and \ref{rmkth2}, on how, following the steps of \cite{CR2}, error estimates could be obtained and, thus Theorems \ref{Teorema-1} and \ref{Teorema-2} could be mathematically proved.

As it is well known, numerical solution of three-dimensional Navier-Stokes equations requires large computational resources, and solving these equations in such a thin domain presents even more numerical problems, while solving the new two-dimensional model presented here is much easier. 

 On the other hand, as we have already mentioned, the new model has the same asymptotic behavior as the Navier-Stokes equations, so, in a certain sense, it encompasses the two limit models presented in subsections \ref{subseccion-4-1} and \ref{subseccion_sw}.
 
For all the above reasons, the new model proposed in this article can be considered a good option for calculating viscous fluid flow between two nearby moving surfaces, without the need to decide a priori whether the flow is typical of a lubrication problem or it is of thin fluid layer type, and without the enormous computational effort that would be required to solve the Navier-Stokes equations in such a thin domain.

We are currently working on performing numerical simulations that allow us to compare the accuracy and computation time required for each of the models that we have mentioned in this article. We hope to be able to present the results of these simulations very soon.

\section*{Acknowledgements}
This work has been partially supported by  the European Union's Horizon 2020 Research and Innovation Programme,
under the Marie Sklodowska-Curie Grant Agreement No 823731 CONMECH. 

\appendix

\section{Coefficients definition} \label{ApendiceA}

In this appendix, we introduce some coefficients that depend either only on the lower bound surface parametrization, $\vec{X}$ or on both the parametrization and the gap $h$. We will use these coefficients throughout this article.

In the first place, the coefficients of the first and second fundamental forms of the surface parametrized by $\vec{X}$ have been denoted by $E, F, G$ and $e,f,g$, respectively:
\begin{eqnarray} E&=&\vec{a}_1 \cdot\vec{a}_1 \label{coef-E} \\
 F&=&\vec{a}_1 \cdot\vec{a}_2 \label{coef-F} \\
 G&=&\vec{a}_2 \cdot\vec{a}_2 \label{coef-G}\\
e&=& -\vec{a}_1 \cdot \dfrac{\partial \vec{a}_{3}}{\partial
\xi_1} = \vec{a}_3 \cdot \dfrac{\partial \vec{a}_{1}}{\partial
\xi_1} \label{coef-e} \\
f&=&-
\vec{a}_1 \cdot \dfrac{\partial \vec{a}_{3}}{\partial
\xi_2}=-\vec{a}_2 \cdot\dfrac{\partial \vec{a}_{3}}{\partial
\xi_1}= \vec{a}_3 \cdot \dfrac{\partial \vec{a}_{1}}{\partial
\xi_2}= \vec{a}_3 \cdot \dfrac{\partial \vec{a}_{2}}{\partial
\xi_1} \label{coef-f} \\
g&=&-\vec{a}_2 \cdot\dfrac{\partial \vec{a}_{3}}{\partial
\xi_2}= \vec{a}_3 \cdot \dfrac{\partial \vec{a}_{2}}{\partial
\xi_2} \label{coef-g} 
\end{eqnarray}
and, from them, we define:

\begin{eqnarray}
A^0&=&\|\vec{a}_1\|^2 \|\vec{a}_2\|^2-\left( \vec{a}_1
\cdot \vec{a}_{2}\right)^2 =EG-F^2=\|\vec{a}_1 \times \vec{a}_2\|^2 \label{A0}  \\
A^1&=&\|\vec{a}_2\|^2 \left(\vec{a}_1 \cdot\dfrac{\partial
\vec{a}_{3}}{\partial \xi_1}\right) + \|\vec{a}_1\|^2
\left(\vec{a}_2 \cdot \dfrac{\partial \vec{a}_{3}}{\partial
\xi_2}\right) \nonumber \\
&&{} - \left( \vec{a}_1 \cdot
\vec{a}_{2}\right)\left(\vec{a}_1 \cdot\dfrac{\partial
\vec{a}_{3}}{\partial \xi_2}+\vec{a}_2 \cdot \dfrac{\partial
\vec{a}_{3}}{\partial \xi_1}\right) = -eG-gE+2fF \label{A^1} \\ 
A^2&=&\left(\vec{a}_1 \cdot\dfrac{\partial
\vec{a}_{3}}{\partial \xi_1} \right)\left(\vec{a}_2 \cdot
\dfrac{\partial \vec{a}_{3}}{\partial \xi_2} \right) -\left(
\vec{a}_1 \cdot\dfrac{\partial \vec{a}_{3}}{\partial \xi_2}\right)
\left( \vec{a}_2 \cdot \dfrac{\partial \vec{a}_{3}}{\partial \xi_1}
\right)=eg-f^2\label{A2}\\
M&=&\begin{pmatrix}G &
-F\\
-F &
E\end{pmatrix}\label{M}
\end{eqnarray}

The following coefficients are involved in the change of variable defined in section \ref{PreviousModels}:
\begin{eqnarray}
\alpha_i
&=&
\alpha_i^0+\varepsilon \xi_3 h \alpha_i^1 +\varepsilon^2 \xi_3^2h^2
\alpha_i^2+\cdots,
\quad (i=1,2) \label{alfaides}\\
\alpha_3 &=&\dfrac{ \xi_3 }{ h}(\alpha_3^0+\varepsilon \xi_3 h \alpha_3^1
+\varepsilon^2 \xi_3^2h^2 \alpha_3^2+\cdots), \label{alfa3des}\\
 \beta_i &=&\beta_i^0+\varepsilon \xi_3 h \beta_i^1 +\varepsilon^2 \xi_3^2h^2
\beta_i^2+\cdots, \quad (i=1,2) \label{betaides}\\
\beta_3&=&\dfrac{ \xi_3 }{ h}(\beta_3^0+\varepsilon \xi_3 h \beta_3^1
+\varepsilon^2 \xi_3^2h^2 \beta_3^2+\cdots), \label{beta3des} \\
\gamma_3&=&\dfrac{1}{\varepsilon h}, \quad \gamma_1 = \gamma_2 = 0, \label{gammades} 
\end{eqnarray}
where 
\begin{eqnarray}
&&\alpha_1^0=\dfrac{\|\vec{a}_2\|^2}{A^0}=\dfrac{G}{EG-F^2}\label{alfa10}\\
&&\alpha_1^1= \dfrac{\vec{a}_2\cdot \dfrac{\partial
\vec{a}_{3}}{\partial \xi_2}- \alpha_1^0 A^1}{ A^0}=- \dfrac{g +
\alpha_1^0 A^1}{ A^0}\label{alfa11}
\\
&&\alpha_1^n= -\dfrac{\alpha_1^{n-2} A^2+\alpha_1^{n-1} A^1}{  A^0},
\quad n\geq 2
\label{alfa1n}\\
&&\alpha_2^0=\beta_1^0=-\dfrac{ \vec{a}_2\cdot \vec{a}_{1}}{ A^0} =-\dfrac{ F}{ A^0}  \label{alfa20_beta10}\\
&&\alpha_2^1=\beta_1^1 =- \dfrac{\vec{a}_2\cdot \dfrac{\partial
\vec{a}_{3}}{\partial \xi_1} +\alpha_2^0 A^1}{ A^0}=
\dfrac{f- \alpha_2^0 A^1}{ A^0} \label{alfa21_beta11}\\
&&\alpha_2^n = \beta_1^n= -\dfrac{ \alpha_2^{n-2} A^2 +
\alpha_2^{n-1} A^1}{ A^0},
\quad n\geq 2  \label{alfa2n_beta1n}\\
&&\alpha_3^0=\dfrac{\dfrac{\partial h}{\partial \xi_2} \vec{a}_{1}
\cdot \vec{a}_2 - \dfrac{\partial h}{\partial \xi_1} \|\vec{a}_2\|^2
}{A^0} =-\alpha_1^0 \dfrac{\partial h}{\partial \xi_1}  -
\alpha_2^0\dfrac{\partial h}{\partial \xi_2}    \label{alfa30}\\
&&\alpha_3^1= \dfrac{\vec{a}_2\cdot\left[\dfrac{\partial h}{\partial
\xi_2}
  \dfrac{\partial
\vec{a}_{3}}{\partial \xi_1} - \dfrac{\partial h}{\partial \xi_1}
\dfrac{\partial \vec{a}_{3}}{\partial \xi_2}\right]-\alpha_3 ^0 A^1
} { A^0} =-\alpha_1^1 \dfrac{\partial h}{\partial \xi_1}  -
\alpha_2^1\dfrac{\partial h}{\partial \xi_2}
\label{alfa31}\\
 && \alpha_3^n=-\dfrac{\alpha_3^{n-2} A^2 + \alpha_3^{n-1} A^1} {A^0},  \quad n\geq 2 \label{alfa3n}
\\
&&\beta_2^0=\dfrac{ \|\vec{a}_{1}\|^2}{ A^0} =\dfrac{ E}{ A^0}  \label{beta20}\\
&&\beta_2^1 = \dfrac{\vec{a}_1\cdot \dfrac{\partial
\vec{a}_{3}}{\partial \xi_1} -\beta_2^0 A^1}{ A^0} =-
\dfrac{ e +\beta_2^0 A^1}{ A^0}\label{beta21}\\
&&\beta_2^n =-\dfrac{ \beta_2^{n-2} A^2 + \beta_2^{n-1} A^1}{ A^0},  \quad n\geq 2   \label{beta22}\\
&&\beta_3^0=\dfrac{\dfrac{\partial h}{\partial \xi_1} \vec{a}_{1}
\cdot \vec{a}_2 - \dfrac{\partial h}{\partial \xi_2}\|\vec{a}_1\|^2
}{A^0} =-\beta_1^0 \dfrac{\partial h}{\partial \xi_1}  -
\beta_2^0\dfrac{\partial h}{\partial \xi_2} \label{beta30}\\
&&\beta_3^1= \dfrac{\dfrac{\partial h}{\partial \xi_1} \left(
\vec{a}_1 \cdot\dfrac{\partial \vec{a}_{3}}{\partial \xi_2}\right) -
\dfrac{\partial h}{\partial \xi_2} \left(  \vec{a}_1
\cdot\dfrac{\partial \vec{a}_{3}}{\partial \xi_1}\right)-\beta_3 ^0
A^1 } { A^0} = -\beta_1^1 \dfrac{\partial h}{\partial \xi_1}  -
\beta_2^1\dfrac{\partial h}{\partial \xi_2}
\label{beta31}\\
 && \beta_3^n=-\dfrac{\beta_3^{n-2} A^2 + \beta_3^{n-1} A^1} {A^0},  \quad n\geq 2
 \label{beta3n}
\end{eqnarray}

The next set of coefficients depend on the parametrization $\vec{X}$ 
and, some of them also depend on function $h$:

\begin{eqnarray}
B^j_{lk}&=& \alpha^j_l (\vec{a}_1 \cdot \vec{a}_k) +
 \beta^j_l (\vec{a}_2 \cdot \vec{a}_k), \quad (j=0,1,2;\quad l=1,2,3;\quad k=1,2, 3) \label{B}\\
C^0_l&=& \alpha_l^0 \left( \vec{a}_1 \cdot \dfrac{\partial \vec{X}}{\partial t}\right) + \beta_l^0 \left( \vec{a}_2 \cdot \dfrac{\partial \vec{X}}{\partial t}\right) \quad (l=1,2,3) \label{C}\\
C^{i,j}_l&=& \alpha_l^i \left( \vec{a}_1 \cdot \dfrac{\partial \vec{X}}{\partial t}\right) + \beta_l^i \left( \vec{a}_2 \cdot \dfrac{\partial \vec{X}}{\partial t}\right) + \alpha_l^j \left( \vec{a}_1 \cdot \dfrac{\partial \vec{a}_3}{\partial t}\right) +
\beta_l^j \left( \vec{a}_2 \cdot \dfrac{\partial \vec{a}_3}{\partial t}\right)\nonumber\\
&& (l=1,2,3; \quad i=1,2; \quad j=0,1,2) \label{Cij}\\
{D^j_{ik}}&=& \alpha_i^j\left(\vec{a}_3\cdot \dfrac{\partial \vec{a}_k}{\partial \xi_1}\right)+\beta_i^j\left(\vec{a}_3 \cdot \dfrac{\partial \vec{a}_k}{\partial \xi_2}\right)\quad (i=1,2,3; \quad k=1,2, 3; \quad j=0,1)\label{D}\\
H^j_{ilk}&=& \alpha_i^j\left(\vec{a}_1\cdot \dfrac{\partial \vec{a}_k}{\partial \xi_l}\right)+\beta_i^j\left(\vec{a}_2 \cdot \dfrac{\partial \vec{a}_k}{\partial \xi_l}\right)\quad (l=1,2; \quad i,k=1,2,3; \quad j=0,1,2)\label{H}\\
I&=&\left(\vec{a}_1 \times \dfrac{\partial
\vec{a}_3}{\partial \xi_2}\right)\cdot \vec{a}_3 + \left(\dfrac{\partial
\vec{a}_3}{\partial \xi_1} \times
\vec{a}_2\right)\cdot \vec{a}_3 \label{I}\\
 J^{i,j}_{lm}&=& \alpha^i_l B^j_{m1} +\beta^i_l B^j_{m2} \quad (l, m=1,2,3; \quad i,j=0,1,2) \label{J}\\
\iota_{lm}^{n}&=& h^{n}\sum_{s=0}^{n}J_{lm}^{s,n-s}, \quad (l,m=1,2; \quad n= 0,1,2,3)\label{iotalm}\\
 {{K}^{j,i}_{l}} &=&\sum_{m=1}^2
 \left(\dfrac{\partial \alpha^j_l}{\partial \xi_m}B^i_{m1}+ \dfrac{\partial \beta^j_l}{\partial \xi_m}  B^i_{m2}+ \alpha^j_l H^i_{mm1} + \beta^j_l H^i_{mm2}\right), \quad (l=1,2,3; \quad i,j=0,1)
 \label{Kji}\\
L_{ikl}^{n,r}&=& h^{n-r} \sum_{m=1}^2 H^0_{imk}
\sum_{s=0}^{n-r} (J_{lm}^{s,n-r-s}+J_{ml}^{s,n-r-s})  \nonumber\\
&&{}+  h^{n-r-1}\delta_{ik} \sum_{s=0}^{n-r} \left[  hK_l^{s,n-r-s} + s \sum_{m=1}^2  \dfrac{\partial h}{\partial \xi_m} J_{lm}^{s,n-r-s}+ 2 r  
J_{3l}^{s,n-r-s}+ s   J_{l3}^{s,n-r-s}\right]\nonumber\\
&&(i,l=1,2;\quad k=1,2,3;\quad  n= 0,1,2,3;\quad  0\leq r\leq n)  \label{Likl_nr}\\
L_{3kl}^{n,r}&=& h^{n-r} \left[ \sum_{m=1}^2 \left(\dfrac{
\partial \vec{a}_{k}  }{ \partial \xi_m} \cdot \vec{a}_3 \right)
\sum_{s=0}^{n-r}( J_{lm}^{s,n-r-s}+J_{ml}^{s,n-r-s} ) +\delta_{k3} \sum_{s=0}^{n-r} K_l^{s,n-r-s}\right]  \nonumber\\
&+&  h^{n-r-1}\delta_{3k} \left[\sum_{s=1}^{n-r} s \dfrac{\partial h}{\partial \xi_m} J_{lm}^{s,n-r-s}+ 2r \sum_{s=0}^{n-r} 
J_{3l}^{s,n-r-s}+\sum_{s=1}^{n-r} s   J_{l3}^{s,n-r-s}\right] \nonumber \\
&& (k=1,2,3;\quad l=1,2; \quad n\geq 0;\quad  0\leq r \leq n) \label{L3kl_nr}\\
\bar{L}_{ikl}^{0,0}& =& L_{ikl}^{0,0}  + \dfrac{1}{h}\dfrac{\partial h}{\partial \xi_k} J^{0,0}_{il} -  \dfrac{ \delta_{ik}}{ h}  J^{0,0}_{3l}, \quad (i,l=1,2;\quad k=1,2,3) \label{Lbar}\\
Q^0_{ik}&=&\alpha_i^0 \left( \vec{a}_1 \cdot \dfrac{\partial \vec{a}_k}{\partial t}\right) + \beta_i^0 \left( \vec{a}_2 \cdot \dfrac{\partial \vec{a}_k}{\partial t}\right)-\sum_{l=1}^2 H^0_{ilk}C^0_l, \quad (i=1,2; \quad k=1,2,3) \label{Q} \\
{Q^0_{3k}}&=&
\left(\vec{a}_3 \cdot \dfrac{
\partial \vec{a}_{k}}{\partial t}\right) -\sum_{l=1}^2
 \left(\vec{a}_3 \cdot \dfrac{
\partial \vec{a}_{k}}{\partial \xi_l}\right) C^0_l \quad (k=1,2)\label{Q3}\\
R^0_{ik} &=&  Q^0_{ik}+ H^0_{ik3} \left( \dfrac{\partial \vec{X}}{\partial t} \cdot \vec{a}_3\right) \quad (i=1,2; \quad k=1,2) \label{R} 
\\
S_{ik}^{n,r}&=&\sum_{s=0}^{n-r}\sum_{l=1}^2 \left[ h^{n-r}\left( \sum_{m=1}^2 \dfrac{
\partial^2   \bar{a}_k}{\partial \xi_l \partial \xi_m} \cdot (\vec{a}_1 \alpha_i^0+\vec{a}_2 \beta_i^0 ) 
 J_{lm}^{s,n-r-s} +  H^0_{ilk}  K_l^{s,n-r-s}\right) \right. \nonumber\\
&+&\left.{}   h^{n-r-1}H^0_{ilk}  \left( 2 r   
J_{3l}^{s,n-r-s} +  s\sum_{m=1}^2 \dfrac{\partial h}{\partial \xi_m} J_{lm}^{s,n-r-s} +s   J_{l3}^{s,n-r-s} \right)\right] \nonumber\\
&+&
 r  \delta_{ik} h^{n-r}\sum_{s=0}^{n-r} \left[ \sum_{m=1}^2 H^{n-r+1}_{mm3}+\dfrac{1}{h} 
 K_3^{s,n-r-s}\right. \nonumber\\
&+&\left.{} 
\dfrac{1}{h^2}\left(
 (s-1) \sum_{m=1}^2 \dfrac{\partial h}{\partial \xi_m}  J_{3m}^{s,n-r-s}+(r+s)   J_{33}^{s,n-r-s}\right)  \right]\nonumber\\
&&(i=1,2,\quad k=1,2,3,\quad n= 0,1,2,3,\quad 0\leq r \leq n ) \label{Sik_nr} \\
S_{3k}^{n,r}&=& h^{n-r} \sum_{s=0}^{n-r}\left[  \sum_{m=1}^2
\left( \sum_{l=1}^2\dfrac{
\partial^2 \vec{a}_{k} }{\partial \xi_l \partial \xi_m} \cdot \vec{a}_3 
 J_{lm}^{s,n-r-s}  \right.\right.\nonumber\\
&+&\left.\left.{}  \left(\dfrac{
\partial \vec{a}_{k} }{ \partial \xi_m} \cdot \vec{a}_3 \right) \left( \dfrac{2  r}{h} 
J_{3m}^{s,n-r-s}+
 K_m^{s,n-r-s} +     
\dfrac{ s }{h}  J_{m3}^{s,n-r-s}+   \dfrac{ s}{h} \sum_{l=1}^2\dfrac{\partial h}{\partial \xi_l} J_{ml}^{s,n-r-s}\right)\right)\right.    \nonumber\\
&+&\left.{}
\delta_{3k} \dfrac{r}{h}\left( K_3^{s,n-r-s}
+  \dfrac{r+s}{h}  
J_{33}^{s,n-r-s}
+ 
\dfrac{s-1}{h} \sum_{m=1}^2 \dfrac{\partial h}{\partial \xi_m}  J_{3m}^{s,n-r-s} +   h 
 \sum_{m=1}^2 H^{n-r+1}_{mm3}\right)
\right]\nonumber\\
&& (k=1,2,3,\quad n\geq 0,\quad 0\leq r \leq n) \label{S3k_nr}\\
\bar{S}_{ik}^{0,0}&=&S_{ik}^{0,0} + \left(\dfrac{  I}{\sqrt{A^0}} -\dfrac{A^1}{A^0}   \right) D^0_{ik} -\dfrac{ 1}{ h} 
\sum_{l=1}^2 J^{0,0}_{3l} H^0_{ilk} - \dfrac{1}{ h \sqrt{A^0}}  \sum_{j=1}^2 J^{0,0}_{ij} \left( \dfrac{\partial  \vec{a}_{k}}{\partial
\xi_j} \cdot \vec{\eta}(h)\right) \nonumber \\
&&\quad (i=1,2; \quad k=1,2,3) \label{Sbar} 
 \\
\vec{\eta}(h)&=&
\dfrac{\partial h}{\partial \xi_2} (\vec{a}_1 \times
\vec{a}_3) + h \left(\vec{a}_1 \times \dfrac{\partial
\vec{a}_3}{\partial \xi_2} + \dfrac{\partial
\vec{a}_3}{\partial \xi_1} \times
\vec{a}_2\right) + \dfrac{\partial h}{\partial
\xi_1} (\vec{a}_3\times \vec{a}_2)\label{eta}\\
\kappa^0_i&=&  \sum_{l=1}^2 \dfrac{
\partial }{ \partial \xi_l}\left( \dfrac{\partial \vec{X}}{\partial t} \cdot \vec{a}_3\right) \left(L^{0,0}_{i3l}-  \dfrac{A^1}{A^0}
 J^{0,0}_{il}\right)
-\dfrac{J^{0,0}_{3i}}{h^2}
   \dfrac{\partial h}{\partial t}  - \dfrac{3}{h}\sum_{k=1}^2 J^{0,0}_{ki} \dfrac{\partial^2 h}{\partial
\xi_k \partial t}\nonumber\\
&+&{}
\left( \dfrac{\partial X}{\partial t} \cdot \vec{a}_3\right) \left[{S}^{0,0}_{i3} -\dfrac{1}{h}  \sum_{l=1}^2 J^{0,0}_{3l}
H^0_{il3} 
 -\dfrac{1}{h \sqrt{A^0} }\sum_{l=1}^2J^{0,0}_{li}\left(  \dfrac{\partial \vec{a}_{3}}{\partial
\xi_l}\cdot 
\vec{\eta}(h)\right)  \right] \quad (i=1,2)\label{kappa} 
\end{eqnarray} 
where $\delta_{ij}$ is the Kronecker Delta.

Finally, we have the coefficients that include the external density of 
volume forces and the friction force

\begin{eqnarray}
F^0_i&=&\int_0^1 f^0_i\, d\xi_3 +\dfrac{s_0}{\rho_0 h}(\vec{f}^1_{R_1} + \vec{f}^1_{R_0}) \cdot \left( \alpha^0_i \vec{a}_1 + \beta^0_i \vec{a}_2 \right) \quad (i=1,2) \label{F} \\
\bar{F}^0_i&=&\dfrac{s_0}{\rho_0 h} \sum_{j=1}^2 J^{0,0}_{ij}\left( \vec{f}^1_{R_1}+ \vec{f}^1_{R_0}\right) \cdot
\vec{a}_j    +\dfrac{  1}{2}   \bar{f}_i^{1,0}
  + {\bar{f}_i^{0,0}} \quad (i=1,2)\label{Fbar2}
\end{eqnarray}

It can be proved that
\begin{equation}
\bar{f}_i^{n,0} = 0 \quad (i=1,2; \ n\ge 1) \label{f_bar_n_0}
\end{equation}
so we obtain 
\begin{equation}
F^0_i = \bar{F}^0_i = \dfrac{s_0}{\rho_0 h} \sum_{j=1}^2 J^{0,0}_{ij}\left( \vec{f}^1_{R_1}+ \vec{f}^1_{R_0}\right) \cdot
\vec{a}_j   + {\bar{f}_i^{0,0}} \quad (i=1,2) \label{Fbar3}
\end{equation}


\begin{thebibliography}{99}

\bibitem{RodTabJMAA2021} J.M. Rodr\'{\i}guez, R. Taboada-V\'azquez, Asymptotic analysis of a thin fluid layer flow between two moving surfaces,  J. Math. Anal. Appl. 507 (1) (2022) 125735, https://doi.org/10.1016/j.jmaa.2021.125735.

\bibitem{Dean1} W.R. Dean, Note on the motion of fluid in a curved pipe, Philos. Mag. 4 (20) (1927) 208–223, https://doi.org/10.1080/14786440708564324.

\bibitem{Dean2} W.R. Dean, The stream-line motion of fluid in a curved pipe (second paper), Philos. Mag. 5 (30) (1928) 673–695, https://doi.org/10.1080/14786440408564513.

\bibitem{FriedrichsDressler} K.O. Friedrichs, R.F. Dressler, A boundary-layer theory for elastic plates, Commun. Pure Appl. Math. 14 (1961) 1–33,
https://doi.org/10.1002/cpa.3160140102.

\bibitem{Goldenveizer} A.L. Goldenveizer, Derivation of an approximated theory of bending of a plate by the method of asymptotic integration of the equations of the theory of elasticity, Prikl. Mat. Meh. 26 (4) (1962) 668–686, https://doi.org/10.1016/0021-8928(62)90161-2.

\bibitem{Rigolot1972} A. Rigolot, Sur une théorie asymptotique des poutres, J. Méc. 11 (4) (1972) 673–703, https://zbmath.org/?q=an:0257. 73013.


\bibitem{CiarletDestuynder1} P.G. Ciarlet, P. Destuynder, A justification of the two dimensional linear plate model, J. Mech. 18 (1979) 315–344, https://zbmath.org/0415.73072.

\bibitem{CiarletDestuynder2} P.G. Ciarlet, P. Destuynder, A justification of a nonlinear model in plate theory, Comput. Methods Appl. Mech. Eng.
17–18 (1979) 227–258, https://doi.org/10.1016/0045-7825(79)90089-6.

\bibitem{Ciarlet1980} P.G. Ciarlet, A justification of the von Kármán equations, Arch. Ration. Mech. Anal. 73 (1980) 349–389, https://doi.org/10.1007/BF00247674.

\bibitem{BermudezViano} A. Berm\'udez, J.M. Via\~no, Une justification des \'equations de la thermo\'elasticit\'e des poutres \`a section variable par des m\'ethodes asymptotiques, RAIRO. Anal. Num\'er. 18 (4) (1984) 347–376, https://doi.org/10.1051/m2an/1984180403471.

\bibitem{TutekAganovicNedelec} Z. Tutek, I. Aganovič, J.-C. Nedelec, A justification of the one-dimensional model of an elastic beam, Math. Methods Appl. Sci. 8 (1986) 502–515, https://doi.org/10.1002/mma.1670080133.

\bibitem{Cimatti1983} G. Cimatti, How the Reynolds equation is related to the Stokes equations, Appl. Math. Optim. 10 (1983) 267–274, https://doi.org/10.1007/BF01448389.

\bibitem{BayadaChambat1986} G. Bayada, M. Chambat, The transition between the Stokes equations and the Reynolds equation: a mathematical proof, Appl. Math. Optim. 14 (1986) 73–93, https://doi.org/10.1007/BF01442229.

\bibitem{ChipotLuskin} G. Chipot, M. Luskin, Existence and uniqueness of solutions to the compressible Reynolds lubrication equation, SIAM J.
Math. Anal. 17 (6) (1986) 1390–1399, https://doi.org/10.1137/0517098.

\bibitem{Cimatti2} G. Cimatti, Existence and uniqueness for non linear Reynolds equation, Int. J. Eng. Sci. 24 (5) (1986) 827–834, https://doi.org/10.1016/0020-7225(86)90116-3.

\bibitem{Cimatti3} G. Cimatti, A rigorous justification of the Reynolds equation, Q. Appl. Math. 45 (1987) 627–644, https://doi.org/10.1090/qam/917014.

\bibitem{Chipot} M. Chipot, On the Reynolds lubrication equation, Nonlinear Anal., Theory Methods Appl. 12 (7) (1988) 699–718, https://doi.org/10.1016/0362-546X(88)90023-5.

\bibitem{Nazarov} S.A. Nazarov, Asymptotic solution of the Navier-Stokes problem on the flow of a thin layer of fluid, Sib. Math. J. 31 (2) (1990) 296–307, https://doi.org/10.1007/BF00970660.

\bibitem{MoiseTemamZiane} I. Moise, R. Temam, M. Ziane, Asymptotic analysis of the Navier-Stokes equations in thin domains, Topol. Methods
Nonlinear Anal. 10 (2) (1997) 249–282, https://projecteuclid.org/euclid.tmna/1476842206.

\bibitem{Grenier} E. Grenier, On the derivation of homogeneous hydrostatic equations, ESAIM: Math. Model. Numer. Anal. 33 (5) (1999)
965–970, \url{http://www.numdam.org/item/M2AN_1999__33_5_965_0/}.

\bibitem{BayadaChambatCiuperca} G. Bayada, M. Chambat, I. Ciuperca, Asymptotic Navier–Stokes equations in a thin moving boundary domain, Asymptot. Anal. 21 (2) (1999) 117–132, https://content.iospress.com/articles/asymptotic-analysis/asy362.

\bibitem{AzeradGuillen} P. Az\'erad, F. Guill\'en, Mathematical justification of the hydrostatic approximation in the primitive equations of geophysical fluid dynamics, SIAM J. Math. Anal. 33 (4) (2001) 847–859, https://doi.org/10.1137/S0036141000375962.

\bibitem{GerbeauPerthame} J.-F. Gerbeau, B. Perthame, Derivation of viscous Saint-Venant system for laminar shallow water; numerical validation,
Discrete Contin. Dyn. Syst., Ser. B 1 (1) (2001) 89–102, https://doi.org/10.3934/dcdsb.2001.1.89.

\bibitem{MarusicPaloka} E. Marušić-Paloka, The effects of flexion and torsion on a fluid flow through a curved pipe, Appl. Math. Optim. 44 (3) (2001) 245–272, https://doi.org/10.1007/s00245-001-0021-y.

\bibitem{Hu} Changbing Hu, Asymptotic analysis of the primitive equations under the small depth assumption, Nonlinear Anal., Theory
Methods Appl. 61 (3) (2005) 425–460, https://doi.org/10.1016/j.na.2004.12.005.

\bibitem{Marche} F. Marche, Derivation of a new two-dimensional viscous shallow water model with varying topography, bottom friction
and capillary effects, Eur. J. Mech. B, Fluids 26 (1) (2007) 49–63, https://doi.org/10.1016/j.euromechflu.2006.04.007.

\bibitem{MarusicPalokaPazanin} E. Marušić-Paloka, I. Pažanin, Fluid flow through a helical pipe, Z. Angew. Math. Phys. 58 (1) (2007) 81–99, https://doi.org/10.1007/s00033-006-0073-6.

\bibitem{DecoeneBMS} A. Decoene, L. Bonaventura, E. Miglio, F. Saleri, Asymptotic derivation of the section-averaged shallow water equations for natural river hydraulics, Math. Models Methods Appl. Sci. 19 (3) (2009) 387–417, https://doi.org/10.1142/S0218202509003474.

\bibitem{PanasenkoStavre} G. Panasenko, R. Stavre, Asymptotic analysis of the Stokes flow in a cylindrical elastic tube, Appl. Anal. 91 (11) (2012) 1999–2027, https://doi.org/10.1080/00036811.2011.584187.

\bibitem{PanasenkoPileckas1} G. Panasenko, K. Pileckas, Asymptotic analysis of the non-steady Navier-Stokes equations in a tube structure. I. The case without boundary layer-in-time, Nonlinear Anal., Theory Methods Appl. 122 (2015) 125–168, https://doi.org/10.1016/j. na.2015.03.008.

\bibitem{PanasenkoPileckas2} G. Panasenko, K. Pileckas, Asymptotic analysis of the non-steady Navier-Stokes equations in a tube structure. II. General case, Nonlinear Anal. 125 (2015) 582–607, https://doi.org/10.1016/j.na.2015.05.018.

\bibitem{SG1} M. Anguiano, F.J. Suárez-Grau, Nonlinear Reynolds equations for non-Newtonian thin-film
fluid flows over a rough boundary, IMA Journal of Applied Mathematics, 84 (1) (2019) 63-95,
https://doi.org/10.1093/imamat/hxy052

\bibitem{SG2} F.J. Suárez-Grau, Asymptotic behavior of a non-Newtonian flow in a thin domain with Navier
law on a rough boundary, Nonlinear Analysis, Theory, Methods and Applications, 117 (2015)
99-123,
https://doi.org/10.1016/j.na.2015.01.013.

\bibitem{SG3} I. Pazanin, F.J. Suárez-Grau, Homogenization of the Darcy-Lapwood-Brinkman flow through
a thin domain with highly oscillating boundaries, Bull. Malays. Math. Sci. Soc., 42 (2019) 3073-3109,
 https://doi.org/10.1007/s40840-018-0649-2.

\bibitem{RTV1} J.M. Rodríguez, R. Taboada-Vázquez, From Navier-Stokes equations to shallow waters with viscosity by asymptotic analysis, Asymptot. Anal. 43 (4) (2005) 267–285, https://content.iospress.com/articles/asymptotic-analysis/asy691.

\bibitem{RTV2} J.M. Rodríguez, R. Taboada-Vázquez, From Euler and Navier-Stokes equations to shallow waters by asymptotic analysis,
Adv. Eng. Softw. 38 (6) (2007) 399–409, https://doi.org/10.1016/j.advengsoft.2006.09.011.

\bibitem{RTV3} J.M. Rodríguez, R. Taboada-Vázquez, A new shallow water model with polynomial dependence on depth, Math. Methods
Appl. Sci. 31 (5) (2008) 529–549, https://doi.org/10.1002/mma.924.

\bibitem{RTV4} J.M. Rodríguez, R. Taboada-Vázquez, A new shallow water model with linear dependence on depth, Math. Comput.
Model. 48 (3–4) (2008) 634–655, https://doi.org/10.1016/j.mcm.2007.11.002.

\bibitem{RTV5} J.M. Rodríguez, R. Taboada-Vázquez, Bidimensional shallow water model with polynomial dependence on depth through
vorticity, J. Math. Anal. Appl. 359 (2) (2009) 556–569, https://doi.org/10.1016/j.jmaa.2009.06.003.

\bibitem{RTV6} J.M. Rodríguez, R. Taboada-Vázquez, Derivation of a new asymptotic viscous shallow water model with dependence on
depth, Appl. Math. Comput. 219 (7) (2012) 3292–3307, https://doi.org/10.1016/j.amc.2011.08.053.































\bibitem{CR1} G. Castiñeira, J.M. Rodríguez, Asymptotic analysis of a viscous fluid in a curved pipe with elastic walls, in: F. Ortegón
Gallego, M. Redondo Neble, J. Rodríguez Galván (Eds.), Trends in Differential Equations and Applications, in: SEMA
SIMAI Springer Series, vol. 8, Springer, Cham, 2016, pp. 73–87, \url{https://doi.org/10.1007/978-3-319-32013-7_5}.

\bibitem{CR2} G. Castiñeira, E. Marušić-Paloka, I. Pažanin, J.M. Rodríguez, Rigorous justification of the asymptotic model describing a curved-pipe flow in a time-dependent domain, Z. Angew. Math. Mech. 99 (1) (2019) 99:e201800154, https://doi.org/10.
1002/zamm.201800154.








\bibitem{CiarletLodsI} P.G. Ciarlet, V. Lods, Asymptotic Analysis of Linearly Elastic Shells. I. Justification of Membrane Shell Equations, Arch. Ration. Mech. Anal. 136 (1996) 119-161, https://doi.org/10.1007/BF02316975.

\bibitem{CiarletLodsII} P.G. Ciarlet, V. Lods, B. Miara, Asymptotic Analysis of Linearly Elastic Shells. II. Justification of Flexural Shell Equations, Arch. Ration. Mech. Anal. 136 (1996) 163–190, https://doi.org/10.1007/BF02316976.

\bibitem{CiarletLodsIV} P.G. Ciarlet, V. Lods, Asymptotic Analysis of Linearly Elastic Shells: Generalized Membrane Shells, Journal of Elasticity 43 (1996) 147-188, https://doi.org/10.1007/BF00042508.

\bibitem{CiarletLodsIII} P.G. Ciarlet, V. Lods, Asymptotic Analysis of Linearly Elastic Shells. III. Justification of Koiter's Shell Equations, Arch. Ration. Mech. Anal. 136 (1996) 191-200, https://doi.org/10.1007/BF02316977.


\bibitem{Weiyan} T. Weiyan, Shallow water hydrodynamics. Elsevier, Amsterdam, 1992.



 
\end{thebibliography}
\end{document}